\documentclass[12pt,psamsfonts,leqno,oneside,letterpaper]{amsart}
\usepackage[dvips,text={6.5truein,9truein},left=1truein,top=1truein]{geometry}
\usepackage{amssymb,amsmath,amscd,graphicx,enumerate}
\usepackage{url}
\let\burl=\url
\usepackage[colorlinks,linkcolor=blue,citecolor=blue,pdfstartview=FitH]{hyperref}
\input xy
\xyoption{all}
\SelectTips{cm}{10}

\usepackage{color}

\theoremstyle{definition}
\newtheorem{para}{}[section]

\newtheorem{remark}[para]{Remark}
\newtheorem{remarks}[para]{Remarks}
\newtheorem{notation}[para]{Notation}
\newtheorem{convention}[para]{Convention}
\newtheorem{definition}[para]{Definition}
\newtheorem{definitions}[para]{Definitions}

\theoremstyle{plain}
\newtheorem{theorem}[para]{Theorem}
\newtheorem{lemma}[para]{Lemma}
\newtheorem{proposition}[para]{Proposition}
\newtheorem{corollary}[para]{Corollary}

\numberwithin{equation}{para}
\numberwithin{figure}{section}

\newcommand\Number{\begin{para}}
\newcommand\EndNumber{\end{para}}
\newcommand\Definition{\begin{definition}}
\newcommand\EndDefinition{\end{definition}}
\newcommand\Definitions{\begin{definitions}}
\newcommand\EndDefinitions{\end{definitions}}
\newcommand\Theorem{\begin{theorem}}
\newcommand\EndTheorem{\end{theorem}}
\newcommand\Remark{\begin{remark}}
\newcommand\EndRemark{\end{remark}}
\newcommand\Remarks{\begin{remarks}}
\newcommand\EndRemarks{\end{remarks}}
\newcommand\Convention{\begin{convention}}
\newcommand\EndConvention{\end{convention}}
\newcommand\Notation{\begin{notation}}
\newcommand\EndNotation{\end{notation}}
\newcommand\Lemma{\begin{lemma}}
\newcommand\EndLemma{\end{lemma}}
\newcommand\Proposition{\begin{proposition}}
\newcommand\EndProposition{\end{proposition}}
\newcommand\Corollary{\begin{corollary}}
\newcommand\EndCorollary{\end{corollary}}
\newcommand\Proof{\begin{proof}}
\newcommand\EndProof{\end{proof}}
\newcommand\Equation{\begin{equation}}
\newcommand\EndEquation{\end{equation}}

\newcommand\Conclusions{\begin{enumerate}[(1)]}
\newcommand\EndConclusions{\end{enumerate}}
\newcommand\Properties{\begin{enumerate}[(1)]}
\newcommand\EndProperties{\end{enumerate}}
\newcommand\Conditions{\begin{enumerate}[(i)]}
\newcommand\EndConditions{\end{enumerate}}

\newcommand\Bullets{\begin{itemize}}
\newcommand\EndBullets{\end{itemize}}

\newdimen\caseindent
\global\caseindent=60pt
\newcommand\Cases{\begingroup}
\newcommand\Case[1]{\par\noindent\hangindent\caseindent
     \hbox to \caseindent{\hskip .5\caseindent minus .5\caseindent
     #1\enspace\hfill}\ignorespaces}
\newcommand\EndCases{\par\endgroup}

\newcommand\haitch{{\mathbb H}}

\newcommand\rk{\mathop{{\rm rk}_2}}

\newcommand\vol{\mathop{\rm Vol}}

\newcommand\Isom{\mathop{\rm Isom}}

\newcommand\genus{\mathop{\rm genus}}
\newcommand\image{\mathop{\rm im}}
\newcommand\Z{{\mathbb Z}}
\newcommand\R{{\mathbb R}}

\newcommand\cala{{\mathcal A}}

\newcommand\calb{{\mathcal B}}
\newcommand\calc{{\mathcal C}}
\newcommand\cald{{\mathcal D}}
\newcommand\calh{{\mathcal H}}

\newcommand\calp{{\mathcal P}}
\newcommand\calq{{\mathcal Q}}
\newcommand\cals{{\mathcal S}}
\newcommand\calt{{\mathcal T}}
\newcommand\calf{{\mathcal F}}
\newcommand\calu{{\mathcal U}}
\newcommand\calv{{\mathcal V}}

\newcommand\calw{{\mathcal W}}
\newcommand\calx{{\mathcal X}}
\newcommand\caly{{\mathcal Y}}

\newcommand\inter{\mathop{\rm int}}
\newcommand\rank{\mathop{\rm rank}}
\newcommand\dist{\mathop{\rm dist}}

\newcommand\DS{{\rm DS}}
\newcommand\ADS{{\rm ADS}}
\newcommand\tg{{\rm TG}}

\newcommand\barchi{\bar\chi}
\newcommand\chibar{\bar\chi}
\newcommand\chiminus{{\chi_{\_{}}}}

\newcommand\hyph{-}

\newcommand{\tC}{\widetilde C}

\newcommand{\tZ}{\widetilde Z}

\newcommand{\tM}{\widetilde M}
\newcommand{\tN}{\widetilde N}

\newcommand\WHAT{{3.66}}

\begin{document}

\author{Ian Agol}
\address{Department of Mathematics\\
University of California, Berkeley\\
970 Evans Hall\\
Berkeley, CA 94720-3840}
\email{ianagol@gmail.com}
\thanks{Partially supported by NSF {grants  DMS-0204142 and DMS-0504975}}

\author{Marc Culler}
\address{Department of Mathematics, Statistics, and Computer Science (M/C 249)\\
University of Illinois at Chicago\\
851 S. Morgan St.\\
Chicago, IL 60607-7045}
\email{culler@math.uic.edu}
\thanks{Partially supported by NSF {grants  DMS-0204142 and DMS-0504975}}

\author{Peter B. Shalen}
\address{Department of Mathematics, Statistics, and Computer Science (M/C 249)\\
University of Illinois at Chicago\\
851 S. Morgan St.\\
Chicago, IL 60607-7045}
\email{shalen@math.uic.edu}

\title{Singular surfaces, mod 2 homology,  and hyperbolic volume, I}

\begin{abstract}
If $M$ is a simple, closed, orientable $3$-manifold such that
$\pi_1(M)$ contains a genus-$g$ surface group, and if $H_1(M;\Z_2)$
has rank at least $4g-1$, we show that $M$ contains an embedded closed
incompressible surface of genus at most $g$.  As an application we
show that if $M$ is a closed orientable hyperbolic $3$-manifold of
volume at most $3.08$, then the rank of $H_1(M;\Z_2)$ is at most
$6$.
\end{abstract}

\maketitle

\section{Introduction and general conventions}

Let $M$ be any closed, orientable, hyperbolic $3$-manifold.  The
volume of $M$ is known to be an extremely powerful topological
invariant, but its relationship to more
classical topological invariants remains elusive.  The main
geometrical result of this paper, Theorem \ref{geom 11}, asserts that
if $\vol M \le 3.08$ then $H_1(M,\Z_2)$ has rank at most $6$.

The Weeks-Hodgson census of closed hyperbolic $3$-manifolds
\cite{snappea} contains two examples, {\tt m135(-1,3)} and {\tt
m135(1,3)}, for which the volume is $<3.08$ and the rank of the first
homology with $\Z_2$ coefficients is $3$. (They are both of volume
$2.666745\ldots$, and they have integer first homology isomorphic to
$\Z_2 \oplus \Z_2 \oplus \Z_4$ and $\Z_2 \oplus \Z_4 \oplus \Z_4$
respectively.) There are no examples in that census for which the
volume is $<3.08$ and the rank of the first homology with $\Z_2$
coefficients is $\ge4$. Thus there is still a substantial gap between
our results and the known examples.  However, the bound on the rank of
$H_1(M;\Z_2)$ given in this paper seems to be better by orders of
magnitude than what could be readily deduced by previously available
methods.
  
The proof of Theorem \ref{geom 11} relies on a purely topological
result, Theorem \ref{top 11}, which states that if $M$ is a closed
$3$-manifold which is simple (see \ref{simple def}), if $\pi_1(M)$ has
a subgroup isomorphic to a genus-$g$ surface group for a given integer
$g$, and if the rank of $H_1(M;\Z_2)$ is at least $4g-1$, then $M$
contains a connected incompressible closed surface of genus $g$.  This
may be regarded as a partial analogue of Dehn's lemma for
$\pi_1$-injective genus-$g$ surfaces.

Theorem \ref{geom 11} will be proved in Section \ref{geometric
  section} by combining Theorem \ref{top 11} with a number
of deep geometric results. These include the Marden tameness
conjecture, recently established by Agol \cite{agol} and by
Calegari-Gabai \cite{cg}; a co-volume estimate for $3$-tame, $3$-free
Kleinian groups due to Anderson, Canary, Culler and Shalen
\cite[Proposition 8.1]{accs}; and a volume estimate for hyperbolic
Haken manifolds recently proved by Agol, Storm and Thurston
\cite{ast}.

The results of \cite{ast} depend in turn on estimates developed by
Perelman in his work \cite{perelman} on geometrization of
$3$-manifolds.

By refining the methods of this paper one can obtain improvements
of Theorems \ref{top 11} and \ref{geom 11}. In particular, in the case
$g=2$, the lower bound of $7$ for the rank of $H_1(M;\Z_2)$ in the
hypothesis of Theorem \ref{top 11} can be replaced by $6$, and the
upper bound of $6$ in the conclusion of Theorem \ref{geom 11} can be
replaced by $5$. The relevant refinements will be explored
systematically in \cite{second}.

Our strategy for proving Theorem \ref{top 11} is based on the method
of two-sheeted coverings used by Shapiro and Whitehead in their proof
\cite{ShapiroWhitehead} of Dehn's Lemma. (This method was inspired by
Papakyriakopoulos's tower construction \cite{Papa}, and was
systematized by Stallings \cite{Stallingsloop}.) We consider a
$\pi_1$-injective genus-$g$ singular surface in the $3$-manifold $M$,
i.e. a map $\phi :K\to M$, where $K$ is a closed orientable genus-$g$
surface, and $\phi _\sharp$ { is} injective. One can construct a
``tower''
\goodbreak
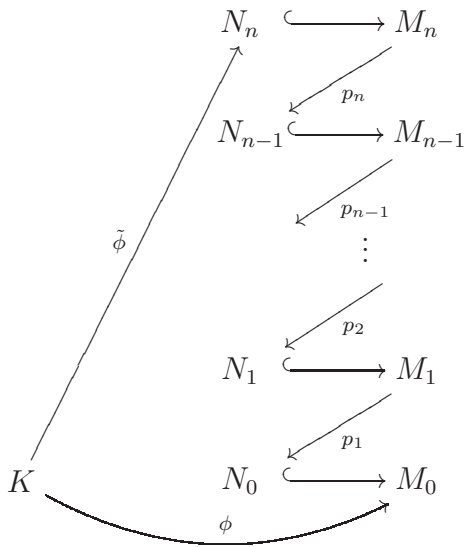
$$
\xymatrix{
&&N_n\hskip8pt\ar@{^{(}->}[r]&M_n\hskip10pt\ar@{->}[dl]^{p_n}\\
&&N_{n-1}\ar@{^{(}->}[r]&M_{n-1}\ar@{->}[dl]^{\genfrac{}{}{0pt}{1}{p_{n-1}}{\vdots}}\\
&&\hbox to 36pt{\strut\hfill}&\hbox to 36pt{\strut\hfill}\ar@{->}[dl]^{p_{2}}\\
&&N_1\hskip8pt\ar@{^{(}->}[r]&M_{1}\hskip10pt\ar@{->}[dl]^{p_{1}}\\
K\ar@/_2pc/[rrr]^\phi\ar@{->}[rruuuu]^{\tilde\phi}&&N_0\hskip8pt\ar@{^{(}->}[r]&M_0\hskip10pt\\
\\}
$$
where the $M_j$ are simple (\ref{simple def}) $3$-manifolds, $N_j$ is
a simple $3$-dimensional submanifold of $M_j$ for $j=0,\ldots,n$, the
$p_j:M_j\to N_{j-1}$ are two-sheeted covering maps,
$\tilde\phi_*:H_1(K;\Z_2)\to H_1(N_n;\Z_2)$ is surjective, and the
diagram commutes up to homotopy. In general this diagram may contain
both closed and bounded manifolds, but we use ideas from
\cite{shalenwagreich} to construct the tower in such a way that if
$H_1(M,\Z_2)$ has rank $\ge 4g-1$, then $H_1(M_j,\Z_2)$ has rank
$\ge 4g-2$  whenever $M_j$ is closed. We also use ideas developed in
\cite{peripheralstructure} based on Simon's results \cite{simon} on
compactification of covering spaces, to construct the tower in such a
way that the (possibly empty and possibly disconnected) surface
$\partial N_j$ is incompressible in $M_j$ for each $j\le n$.

The manifold $N_n$ always has non-empty boundary. This is obvious if
$\partial M_n\ne\emptyset$. If $M_n$ is closed then $H_1(M_n;\Z_2)$ has
rank at least  $4g-2$,  whereas the surjectivity of
$\tilde\phi_*:H_1(K;\Z_2)\to H_1(N_n;\Z_2)$ implies that the rank of
$H_1(N_n;\Z_2)$ is at most  $2g$. It follows that in this case
$N_n$ is a proper submanifold of $M_n$, and hence $\partial
N_n\ne\emptyset$.

We in fact show, using elementary arguments based on
Poincar\'e-Lefschetz duality, that if the map
$\tilde\phi_*:H_2(K;\Z)\to H_2(N_n;\Z)$ is trivial, then $\partial
N_n$ has a component $F$ of genus at most $g$.  In the case where
$\tilde\phi_*:H_2(K;\Z)\to H_2(N_n;\Z)$ is non-trivial, we use Gabai's
results from \cite{gabai} to show that $N_n$ contains a non-separating
incompressible closed surface $F$ of genus at most $g$.

The rest of the proof consists of showing that if a given $M_j$, with
$0<j\le n$, contains a closed incompressible surface of genus at most
$g$, then $N_{j-1}$ also contains such a surface. The surface in
$N_{j-1}$ will be incompressible in $M_{j-1}$, as well as in
$N_{j-1}$, because $\partial N_{j-1}$ is incompressible in
$M_{j-1}$. It is at this step that we need to know that closed
manifolds in the tower have first homology with $\Z_2$-coefficients of
rank at least $4g-2$. Indeed, Proposition \ref{new improved old prop
  3} implies that the existence of a closed incompressible surface of
genus at most $g$ in a $2$-sheeted covering of a simple compact
$3$-manifold $N$ implies the existence of such a surface in $N$ itself
unless $N$ is closed and $H_1(N;\Z_2)$ has rank at most $4g-3$.

Proposition \ref{new improved old prop 3} involves the notion of a
``book of $I$-bundles'' which we define formally in \ref{true book
def}.  Books of $I$-bundles in PL $3$-manifolds arise naturally as
neighborhoods of ``books of surfaces'' (\ref{bosdef}).  We may think
of a book of surfaces as being constructed from a $2$-manifold with
boundary $\hat\Pi$, whose components have Euler characteristic $\le0$,
and a closed $1$-manifold $\Psi$, by attaching $\partial\hat\Pi$ to
$\Psi$ by a covering map. The components of $\Psi$ and
$\Pi=\inter\hat\Pi$ are respectively ``bindings'' and ``pages.''  A
book of $I$-bundles comes equipped with a corresponding decomposition
into ``pages'' which are $I$-bundles over surfaces, and ``bindings''
which are solid tori. (In the informal discussion that we give in this
introduction, the extra structure defined by the decomposition will be
suppressed from the notation.)

With these notions as background we shall now sketch the proof of
Proposition \ref{new improved old prop 3}.  An incompressible surface
$F$ in a two-sheeted covering space of $N$, if it is in general
position, projects to $N$ via a map which has only double-curve
singularities. After routine modifications one obtains a map from $F$
to $N$ with the additional property that its double curves are
homotopically non-trivial. In particular, the image of such a map is a
book of surfaces $X$.  A regular neighborhood $W$ of $X$ in $N$ is
then a book of $I$-bundles, which has Euler characteristic $\ge2-2g$
if $F$ has genus at most $g$. Using the the simplicity of $N$ one can
then produce a book of $I$-bundles $V$ with $W\subset
V\subset N$ and $\chi(W)\ge2-2g$, such that each page of $W$ has
strictly negative Euler characteristic. (This step is handled by
Lemma \ref{make book}.)

We now distinguish two cases. In the case where some page $P_0$ of $V$
has the property that $P_0\cap\partial V$ is contained in a single
component of $\partial V$, we show that by splitting bindings of the
book of surfaces $X$, one can produce an embedded (possibly
disconnected) closed, orientable surface $S$ which is homologically
non-trivial in $N$.  Ambient surgery on $S$ in $N$ then produces a
non-empty incompressible surface whose components have genus at most
$g$. In the case where no such page $P_0$ exists, an Euler
characteristic calculation shows that the boundary components of $V$
have genus at most $g$. In this case, ambient surgery on $\partial V$
produces a non-empty incompressible surface whose components have
genus at most $g$. We show that this surface is non-empty unless $V$
carries $\pi_1(N)$. But for a book of $I$-bundles $V$ whose Euler
characteristic is at least $2-2g$, and whose pages are all of negative
Euler characteristic, one can show that $H_1(V;\Z_2)$ has rank at most
$4g-3$ (this is included in Lemma \ref{moosday}); so in the case where
$V$ carries $\pi_1(N)$, the rank of $H_1(N;\Z_2)$ is at most $4g-3$.

The details and background needed for the proof of Proposition
\ref{new improved old prop 3} occupy Sections \ref{boib section
1}---\ref{down}. Section \ref{singularity section} provides the
combinatorial background needed to construct the tower, while Sections
\ref{homology section} and \ref{new section} provide the homological
background. The application of Gabai's results mentioned above appears
in Section \ref{new section}. The material on towers proper, and the
proof of the main topological theorem and its corollary, are given in
Section \ref{tower section}, and the geometric applications are given
in Section \ref{geometric section}.

The rest of this introduction will be devoted to indicating some
conventions that will be used in the rest of the paper.

\Number
In general, if $X$ and $Y$ are subsets of a set, we denote by
$X\setminus Y$ the set of elements of $X$ that do not belong to
$Y$. In the case where we know that $Y\subset X$ and wish to emphasize
this we will write $X-Y$ for $X\setminus Y$.
\EndNumber

\Number
A {\it manifold} may have a boundary. If $M$ is a manifold, we shall
denote the boundary of $M$ by $\partial M$ and its interior
$M-\partial M$ by $\inter M$.

In many of our results about manifolds of dimension $\le3$ we do not
specify a category. These results may be interpreted in the category
in which manifolds are topological, PL or smooth, and submanifolds are
respectively locally flat, PL or smooth; these three categories are
equivalent in these low dimensions as far as classification is
concerned. In much of the paper the proofs are done in the PL
category, but the applications to hyperbolic manifolds in Section
\ref{geometric section} are carried out in the smooth category.
\EndNumber

\Number\label{separating def}
A (possibly disconnected) codimension-$1$ submanifold $S$ of a
manifold $M$ is said to be {\it separating} if $M$ can be written as
the union of two $3$-dimensional submanifolds $M_1$ and $M_2$ such
that $M_1\cap M_2 = S$.
\EndNumber

\Number\label{abcd goldfish}
We shall say that a map of topological
spaces $f:\calx\to \caly$ is {\it $\pi_1$-injective} if for every path
component $X$ of $\calx$, the map $f|X$ induces an injection from
$\pi_1(X)$ to $\pi_1(Y)$, where $Y$ is the path component of $\caly$
containing $f(X)$.
We shall say that a subset $A$ of a 
topological space $X$ is {\it $\pi_1$-injective} in $X$ if the
inclusion map $A\to X$ is $\pi_1$-injective.
\EndNumber

\Number\label{oiler}
If $X$ is a space having the homotopy type of a finite CW complex, the
Euler characteristic of $X$ will be denoted by $\chi(X)$.  We have
$\chi(X)=\sum_{j\in\Z}\dim_FH_j(X;F)$ for {\it any} field $F$: the sum
is independent of $F$ by virtue of the standard observation that it is
equal to $\sum_{j\in\Z}(-1)^jc_j$, where $c_j$ denotes the number of
$j$-cells in a finite CW complex homotopy equivalent to $X$.

We shall often write $\barchi(X)$ as shorthand for $-\chi(X)$.
\EndNumber

\Number If $x$ is a point of a compact PL space $X$, there exist a
finite simplicial complex $K$ and a PL homeomorphism $h:X\to|K|$ such
that $v=h(x)$ is a vertex of $K$. If $L$ denotes the link of $v$ in
$K$ then the PL homeomorphism type of the space $|L|$ depends only on
$X$ and $x$, not on the choice of $K$ and $h$.  We shall refer to $L$
as the {\it link of $x$ in $X$}, with the understanding that it is
defined only up to PL homeomorphism.  \EndNumber

\Number\label{Pi notation}
Suppose that $x$ is a point of a compact PL space $X$ and that
$n\ge0$ is an integer. The link of $x$ is PL homeomorphic to $S^{n-1}$
if and only if $x$ is an $n$-manifold point of $X$, i.e. some
neighborhood of $x$ is piecewise-linearly homeomorphic to ${\bf R}^n$.

If $X$ is a compact PL space of dimension at most $2$, we
shall denote by $\Pi(X)$ the set of all $2$-manifold points of
$X$. Note that $\Pi(X)$ is an open subset of $X$, and with its induced
PL structure it is a PL $2$-manifold. Furthermore, $X-\Pi(X)$ is a
compact PL subspace of $X$.
\EndNumber

\Number
Let $F$ be a properly embedded orientable surface in an orientable $3$-manifold
$M$. We define a {\it compressing disk} for $F$ in $M$ to be a disk
$D\subset M$ such that $D\cap F=\partial D$, and such that $\partial
D$ is not the boundary of a disk in $F$. It is a standard consequence
of the loop theorem that $F$ is $\pi_1$-injective in $M$ if and only
if there is no compressing disk for $F$ in $M$.

A closed orientable surface $S$ contained in the interior of an
orientable $3$-manifold $M$ will be termed {\it incompressible} if $S$
is $\pi_1$-injective in $M$ and no component of $S$ is a $2$-sphere.
(We have avoided using the term ``incompressible'' for surfaces that
are not closed.) 
\EndNumber

\Number
An {\it essential arc} in a $2$-manifold $F$ is a properly embedded
arc in $F$ which is not the frontier of a disk.
\EndNumber

\Definitions\label{simple def}
A $3$-manifold $M$ will be termed {\it irreducible} if every
$2$-sphere in $M$ bounds a ball in $M$.  We shall say that $M$ is {\it
boundary-irreducible} if $\partial M$ is $\pi_1$-injective in $M$, or
equivalently if, for every properly embedded disk $D\subset M$, the
simple closed curve $\partial D$ bounds a disk in $\partial M$. We
shall say that a $3$-manifold $M$ is {\it simple} if (i) $M$ is
compact, connected, orientable, irreducible and boundary-irreducible;
(ii) no subgroup of $\pi_1(M)$ is isomorphic to $\Z\times\Z$; and
(iii) $M$ is not a closed manifold with finite fundamental
group.
\EndDefinitions

\Number\label{simple covering}
It is a theorem due to Meeks, Simon and Yau \cite{MSY} that a covering
space of an irreducible orientable $3$-manifold is always irreducible.
Given this result, it follows formally from our definition of
simplicity that if a compact, orientable $3$-manifold $M$ is simple,
then every connected finite-sheeted covering of $M$ is also simple.
\EndNumber

\Number\label{horizontal def}
The unit interval $[0,1]$ will often be denoted by $I$.

By an {\it $I$-bundle} we shall always mean a compact space equipped
with a specific locally trivial fibration over some (often unnamed)
base space, in which the fibers are homeomorphic to $[0,1]$.  (The
reader is referred to \cite[Chapter 10]{hempel} for a general
discussion of $3$-dimensional $I$-bundles.)

By a {\it Seifert fibered manifold} we shall always mean a compact
$3$-manifold equipped with a specific Seifert fibration.

In particular, the notion of a {\it fiber} of an $I$-bundle or a
Seifert fibered manifold is well defined, although the fiber
projection and base space will often not be explicitly named.  A
compact subset of an $I$-bundle or Seifert fibered space will be
called {\it horizontal} if it meets each fiber in one point. A compact
set will be called {\it vertical} if it is a union of fibers.
  
If $\calp$ is an $I$-bundle, we define the {\it horizontal boundary}
of $\calp$ to be the subset of $\calp$ consisting of all endpoints of
fibers of $\calp$. We shall denote the horizontal boundary of $\calp$
by $\partial_h\calp$.
  
In the case where the base of the $I$-bundle $\calp$ is a $2$-manifold
$F$ (so that $\calp$ is a $3$-manifold), we define the {\it vertical
boundary} of $\calp$ to be $p^{-1}(\partial F)$, where $p:\calp\to F$
denotes the bundle map.  Note that in this case we have
$\partial\calp=\partial_v\calp\cup\partial_h\calp$, and if $\calp$ is
orientable then $\partial_v\calp$ is always a finite disjoint union of
annuli.
\EndNumber

\Number
The {\it rank} of a finitely generated group $\Gamma$ is the
cardinality of a minimal generating set for $\Gamma$. In particular,
the trivial group has rank $0$.

A group $\Gamma$ is said to be {\it freely indecomposable}
if $\Gamma$ is not the free product of two non-trivial subgroups.
\EndNumber

\Number
If $V$ is a finite dimensional vector space over $\Z_2$ then the
dimension of $V$ will be denoted $\rk V$.  If $X$ is a topological
space, we will set $\rk X = \rk H_1(X;\Z_2)$.
\EndNumber

\section{Books of \texorpdfstring{$I$}{I}-bundles}
\label{boib section 1}

\Definition
A {\it generalized book of $I$-bundles} is a triple
$\calw=(W,\calb,\calp)$, where $W$ is a (possibly empty) compact,
orientable $3$-manifold, and $\calb,\calp\subset W$ are submanifolds
such that
\Bullets
\item $\calb$ is a (possibly disconnected) Seifert fibered space,
\item $\calp$ is an $I$-bundle over a (possibly disconnected)
$2$-manifold, and every component of $\calp$ has Euler characteristic $\le0$,
\item $W=\calb\cup\calp$,
\item $\calb\cap\calp$ is the vertical boundary of $\calp$, and
\item $\calb\cap\calp$ is
vertical in the Seifert fibration of $\calb$.
\EndBullets

We shall denote $W$, $\calb$ and $\calp$ by $|\calw|$, $\calb_\calw $
and $\calp_\calw $ respectively. The components of $\calb_\calw $ will
be called {\it bindings} of $\calw$, and the components of
$\calp_\calw$ will be called its {\it pages}.  The submanifold
$\calp\cap\calb$, whose components are properly embedded annuli in
$W$, will be denoted $\cala_\calw$.

If $B$ is a binding of a generalized book of $I$-bundles $\calw$,
we define the {\it valence} of $B$ to be the number of components of
$\cala_\calw$ that are contained in $\partial B$. 

A generalized book of $I$-bundles $\calw$ will be termed {\it
  connected} if the manifold $|\calw|$ is connected. Likewise, $\calw$
will be termed {\it boundary-irreducible} if $|\calw|$ is
boundary-irreducible.
\EndDefinition

\Definitions\label{true book def}
A {\it book of $I$-bundles} is a generalized book of $I$-bundles
$\calw$ such that
\Bullets
\item $|\calw|\ne\emptyset$,
\item each binding of $\calw$ is a solid torus, and
\item each binding of $\calw$ meets at least one page of $\calw$.
\EndBullets

If $B$ is a binding of a book of $I$-bundles $\calw$, there is a
unique integer $d>0$ such that for every component $A$ of
$\cala_\calw$ contained in $\partial B$, the image of the inclusion
homomorphism $H_1(A;\Z)\to H_1(B;\Z)$ has index $d$ in $H_1(B;\Z)$. We
shall call $d$ the {\it degree} of the binding $B$.
\EndDefinitions

\Lemma\label{no annular pages}
Suppose that $\calw$ is a generalized book of $I$-bundles. Then
there is a generalized book of $I$-bundles $\calw_0$ such that
\Conclusions
\item\label{good morrow} $|\calw_0|=|\calw|$,
\item\label{i prithee discover} every page of $\calw_0$ has strictly
negative Euler characteristic, and
\item\label{steal purchase or borrow}every page of $\calw_0$ is
  a page of $\calw$.
\EndConclusions
\EndLemma

\Proof
Set $W=|\calw|$, $\calb=\calb_\calw$ and $\calp=\calp_\calw$. Let
$\calq$ denote the union of all components $P$ of $\calp$ such that
$\chi(P)=0$.  Then $\calq$ is an $I$-bundle over a compact surface $A$
whose components are annuli and M\"obius bands, and $\calq\cap\calb$
is the induced $I$-bundle over $\partial A$. Hence every component $Q$
of $\calq$ is a solid torus, and $Q\cap\calb$ consists of either a
single annulus of degree $2$ in $Q$, or of two annuli of degree $1$ in
$Q$.  Since every such annulus is also vertical in the Seifert
fibration of $\calb$, it follows that this Seifert fibration may be
extended to a Seifert fibration of the manifold $\calb_0=\calb\cup
\calq$, in such a way that each component of $\calq$ contains either
no singular fiber, or exactly one singular fiber of order
$2$. Furthermore, since every component of $\calq$ meets $\calb$,
every component of $\calb_0$ contains a component of $\calb$.

The manifold $\calp_0=\calp-\calq$ is a union of components of $\calp$
and therefore inherits an $I$-bundle structure. It is now clear that
$\calw_0 = (W,\calb_0,\calp_0)$ is a generalized book of
$I$-bundles. It follows from the definition of $\calq$ that $\calw_0$
satisfies conclusions (\ref{i prithee discover}) and (\ref{steal
purchase or borrow}) of the lemma.   \EndProof

\Lemma\label{Seifert holes}
Suppose that $\hat B$ is a connected, Seifert-fibered submanifold of a
simple, closed, orientable $3$-manifold $M$. Then either
\Conclusions
\item\label{hole one} $\hat B$ is a solid torus, or
\item\label{hole two} $\hat B$ is contained in a ball in $M$, or
\item\label{hole three} some component of $M-\inter \hat B$ is a solid torus.
\EndConclusions
\EndLemma

\Proof Since $M$ is simple and $\hat B$ is Seifert-fibered, we have
$\hat B\ne M$, i.e. $\partial \hat B\ne\emptyset$. Since the
components of $\partial \hat B$ are tori and $M$ is simple, $\partial
\hat B$ cannot be $\pi_1$-injective in $M$. Hence there is a
compressing disk for $\partial \hat B$ in $M$.  If $D\subset \hat B$
then $\hat B$ is a boundary-reducible Seifert fibered space and hence
(\ref{hole one}) holds. The other possibility is that $D\cap \hat
B=\partial D$. In this case, let $V$ denote a regular neighborhood of
$D$ relative to $M-\inter \hat B$. The boundary of the manifold $\hat
B\cup V$ has a unique sphere component $S$. Since $M$ is irreducible,
$S$ bounds a ball $\Delta\subset M$. We must have either
$\Delta\supset \hat B$, which gives conclusion (\ref{hole two}), or
$\inter \Delta\cap \hat B=\emptyset$; in the latter case, $\Delta\cup
V$ is a solid torus component of $M-\inter \hat B$, and so (\ref{hole
  three}) holds.
\EndProof

\Lemma\label{make book} Suppose that $M$ is a simple, closed,
orientable $3$-manifold, and that $\calw$ is a connected generalized
book of $I$-bundles such that $W=|\calw|\subset M$. Suppose that
$\chi(W)<0$, and that $\calp_\calw $ is $\pi_1$-injective in
$M$. Then there is a connected book of $I$-bundles $\calv$
with $V=|\calv|\subset M$, such that
\Conclusions
\item\label{a pair of pizza pies} $V\supset W$,
\item\label{Euler doesn't care}$\chibar(V)=\chibar(W)$, 
\item\label{he does and he doesn't}$\chi(P)<0$ for every page $P$
  of $\calv$,
\item\label{wysiwyg} $\partial V$ is a union of components of
  $\partial W$,
\item\label{little birdies' dirty feet} every component of $\overline{V-W}$ is a solid torus,
\item\label{and me without a spoon} every page  of
  $\calv$ is a page   of  $\calw$, and 
\item\label{mutilated monkey meat}for each page $P$ of $\calv$ we have
  $P\cap\partial V=P\cap\partial W$.
\EndConclusions
\EndLemma

\Proof
Let $\calw_0 = (W, \calb, \calp)$ be a generalized book of
$I$-bundles satisfying conditions (\ref{good morrow})--(\ref
{steal purchase or borrow}) of Lemma \ref{no annular pages}.
Since each page of $\calw_0$ is also a page of $\calw$, the hypothesis
implies that each page of $\calw_0$ is $\pi_1$-injective in $M$.

Let $B$ be any binding of $\calw_0$.  We will show that the
Seifert fibers of $B$ are homotopically non-trivial in $M$.  Since
$\calw_0$ is connected and $\chi(|\calw_0|)<0$, the binding $B$ must
meet some page $P$ of $\calw_0$.  Let $A$ be one of the annulus
components of $B\cap P$.  Then $A$ is a component of the vertical
boundary of $P$ and, since $\chi(P)<0$, it follows that $A$ is
$\pi_1$-injective in $P$.  Since $P$ is $\pi_1$-injective in $M$, it
follows that $A$ is also $\pi_1$-injective in $M$. Recalling that the
annulus $A$ is saturated in the Seifert fibration of $B$, we may
conclude that each Seifert fiber of $B$ is homotopically non-trivial
in $M$.

Now for any binding $B$ of $\calw_0$ let us define $\hat B$ to be the
union of $B$ with all of the solid torus components of
$\overline{M-B}$.  We will show that $\hat B$ is a Seifert fibered
submanifold of $M$ such that $\hat B\cap\calw_0 = B$.  If $J$ is any
solid torus component of $\overline{M-B}$ then no page of $\calv_0$
can be contained in $J$, since the pages are $\pi_1$-injective in $M$
and have negative Euler characteristic.  Thus $\inter J$ must be
disjoint from all of the pages of $\calw_0$.  This implies that $\hat
B\cap\calw_0 = B$.  If $F\subset\partial J$ is a fiber of the Seifert
fibered space $B$ then, since $F$ is homotopically non-trivial in $M$,
the simple closed curve $F\subset\partial J$ cannot be a meridian
curve for the solid torus $J$. It follows that the Seifert fibration
of $B$ may be extended to a Seifert fibration of $B=B\cup J$, and
hence that $\hat B$ admits a Seifert fibration.

Next we will show that $\hat B$ is, in fact, a solid torus.  We know
that $\hat B$ must satisfy one of the conditions (1)---(3) of Lemma
\ref{Seifert holes}.  Condition \ref{Seifert holes}(\ref{hole three})
is ruled out since, by construction, no component of $M-\inter \hat B$
is a solid torus.  The fact that the Seifert fibers of $B$ are
homotopically non-trivial in $M$ implies that the inclusion
homomorphism $\pi_1(B)\to\pi_1(M)$ has non-trivial image and thus $B$
cannot be contained in a ball in $M$.  This rules out condition
\ref{Seifert holes}(\ref{hole two}).  Thus we conclude that condition
\ref{Seifert holes}(\ref{hole one}) must hold, i.e that $\hat B$ is a
solid torus.

Since each binding of $\calw$ must meet some page, and since no page
can be contained in a solid torus, we have that if $B_1$ and $B_2$ are
distinct bindings of $\calw_0$, then $\hat B_1$ is disjoint from $\hat
B_2$.  We define $\calb'$ to be the union of the solid tori $\hat
B$ as $B$ ranges over all bindings of $\calw_0$, and we set $V =
\calb'\cup\calp$. We have $\calb'\cap|\calw_0| = \calb$.  It
follows that $\calv = (V,\calb',\calp)$ is a book of $I$-bundles,
and that every page of $\calv$ has strictly negative Euler
characteristic.

We shall now complete the proof by observing that $V$ satisfies
Conclusions (\ref{a pair of pizza pies})---(\ref{mutilated monkey
  meat}) of the present lemma.  Conclusions (\ref{a pair of pizza
  pies}), (\ref{wysiwyg}) and (\ref{little birdies' dirty feet}) are
immediate from the construction of $V$, and they imply Conclusion
(\ref{Euler doesn't care}).  The pages of $\calv$ are the same as the
pages of $\calw_0$, and each page of $\calw_0$ is a page of $\calw$
and has negative Euler characteristic. Hence Conclusions (\ref{he does
  and he doesn't}) and (\ref{and me without a spoon}) hold.  Since
$\partial W$ is the union of $\partial V$ with a collection of tori
that are disjoint from all pages, it follows that $P\cap\partial V =
P\cap\partial W$ for every page $P$ of $\calv$.  This is Conclusion
(\ref{mutilated monkey meat}).
\EndProof

Recall that in \ref{Pi notation} we defined $\Pi(X)\subset X$ to be
the set of $2$-manifold points in an arbitrary compact PL space $X$ of
dimension at most $2$, and we observed that $X-\Pi(X)$ is a compact PL
subset of $X$. It follows that $\Pi(X)$ has the homotopy type of a
compact PL space. In particular $\chi(\pi)$ is a well-defined integer
for every component $\pi$ of $\Pi(X)$.

\Definition\label{bosdef}
We define a {\it book of surfaces} to be a compact PL space $X$ such
that
\Properties
\item the link of every point of $x\in X$ is PL homeomorphic to
the suspension of some non-empty finite set $Z_x$; and
\item for every component $\pi$ of $\Pi(X)$ we have $\chi(\pi)\le0$.
\EndProperties
The cardinality of the set $Z_x$ appearing in condition (1) is clearly
uniquely determined by the point $x$. It will be called the {\it
order} of $x$.
\EndDefinition

\Number
Note that a point $x$ in a book of surfaces $X$ has order $2$ if and
only if $x\in\Pi(X)$. 

It also follows from the definition that if $X$ is a book of surfaces,
the set $X-\Pi(X)$ is a compact PL $1$-manifold, which will be denoted
by $\Psi(X)$. The components of $\Psi(X)$ and $\Pi(X)$ may be
respectively thought of as {\it bindings} and {\it pages} of $X$.

We also observe that if $M$ is a PL $3$-manifold and if $S_1$ and
$S_2$ are closed surfaces in $\inter M$ which meet transversally, then
$S_1\cup S_2$ is a book of surfaces.
\EndNumber

\Lemma\label{book structure}
If $X$ is a book of surfaces, there exist a (possibly disconnected)
compact PL surface $F$ and a PL map $r:F\to X$ such that
\Conclusions
\item $r|\inter F$ is a homeomorphism of $\inter F$ onto $\Pi(X)$,
and
\item $r|\partial F$ is a covering map from $\partial F$ to $\Psi(X)$.
\EndConclusions
\EndLemma

\Proof
Let us identify $X$ with $|K|$, where $|K|$ is some finite simplicial
complex. After subdividing $K$ if necessary we may assume that for
every closed simplex $\Delta$ of $K$ the set $\Delta\cap\Psi(X)$ is a
(possibly empty) closed face of $\Delta$.
 
Let $\cald$ denote the abstract disjoint union of all the closed
$2$-simplices of $X$, and let $i:\cald\to X$ denote the map which is
the inclusion on each closed $2$-simplex. For each point $z\in\cald$
let $\Delta_z$ denote the closed $2$-simplex containing $z$. We define
a relation $\sim$ on $\cald$ by writing $z\sim w$ if and only if (i)
$\Delta_z\cap\Delta_w\not\subset\Psi(X)$ and (ii) $i(z)=i(w)$.  It is
straightforward to show that $\sim$ is an equivalence relation.  The
quotient space $F=\Delta/\sim$ inherits a PL structure from
$\cald$. The definition of $\sim$ implies that there is a unique map
$r:F\to X$ such that $r\circ q=i$, where $q:\cald\to F$ is the
quotient map, and that $r$ maps $E=r^{-1}\Pi(X)$ homeomorphically onto
$\Pi(X)$.

If $x$ is a point of $\Psi(X)$, then since $X$ is a book of surfaces,
there exist a neighborhood $A$ of $x$ in $\Psi(X)$, and a neighborhood
$V$ of $x$ in $X$, such that $A$ is a PL arc, $V$ is a union of PL
disks $D_1\cup\cdots\cup D_m$, where $m$ is the order of $x$ in $X$,
and $D_i\cap D_j=A$ whenever $i\ne j$. The definition of $\sim$
implies that $r^{-1}(V)$ is a disjoint union of PL disks $\widetilde
D_1,\ldots,\widetilde D_m$ such that $r$ maps $\widetilde D_i$
homeomorphically onto $D_i$ for $i=1,\ldots,m$. Hence $F$ is a PL
surface with interior $E$ and boundary $r^{-1}(\Psi(X))$, and
$r|\partial F:\partial F\to\Psi(X)$ is a covering map.
\EndProof

\Lemma\label{neighborhood is a book} Suppose that $M$ is an orientable
PL $3$-manifold and that $X\subset\inter M$ is a book of
surfaces. Then there is a book of $I$-bundles $\calw$ such that
\Conclusions
\item \label{reg nbhd}$|\calw|=W$ is a regular neighborhood of $X$;
\item \label{more reg nbhd}$|\calb_\calw|$ is a regular neighborhood
  of $\psi(X)$;
\item \label{gopher}for every page $P$ of $\calw$, the set $X\cap
  P$  is a section of the $I$-bundle $P$; and
\item\label{twofer} $\calp_\calw$ is a regular neighborhood
  in $M$ of a deformation retract of  $\Pi(X)$.
\EndConclusions
\EndLemma

\Proof
Let $\calb$ be a regular neighborhood of $\Psi(X)$ in $M$ such that
$N=\calb\cap X$ is a regular neigborhood of $\Psi(X)$ in the PL space
$X$. Every component of $\calb$ is a solid torus. Since $\Pi(X)$ is an
open $2$-manifold, $Y=X\cap\overline{M-\calb}$ is a compact
$2$-manifold and a deformation retract of $\Pi(X)$. In particular, in
view of condition (2) in the definition of a book of surfaces, every
component of $Y$ has Euler characteristic $\le0$. Let $\calp$ be a
regular neighborhood of $Y$ in $\overline{M-N}$.  Then
$W=\calb\cup\calp$ is a regular neighborhood of $X$ in $M$. We may
give $\calp$ the structure of an $I$-bundle over $Y$ in such a way
that $Y$ is identified with a section of the bundle. We have
$\calp\cap\calb=\partial_v\calp$, and $\chi(P)\le0$ for every
component $P$ of $\calp$.

Let $F$ be the surface, and $r:F\to X$ the map, given by Lemma
\ref{book structure}.  We have $N=r(C)$, where $C$ is a
collar neighborhood of $\partial F$ in $F$.  Now if $A$ is any
component of $\partial_v\calb$, then $A\cap Y$ is a component of
$\partial Y$ and therefore cobounds an annulus component of $C$ with
some component $\widetilde \psi_A$ of $\partial F$.  It follows from
\ref{book structure} that $r|\widetilde \psi_A$ is a covering map of some
degree $d_A$ to some component $\psi_A$ of $\psi(X)$. The annulus $A$
lies in the boundary of the component $B_A$ of $\calp$ containing
$\psi_A$, and the (unsigned) degree of $A$ in the solid torus $B_A$ is
$d_A$. In particular, every component of $\partial_v\calp$ has
non-zero degree in the component of $\calb$ containing it. This
implies that $\calw=(W,\calb,\calp)$ is a book of $I$-bundles.

Each page $P$ of $\calp$ was constructed as an $I$-bundle over a
component $Y_0$ of $Y$, where $Y_0$ is identified with a section of
the bundle.  Since $Y_0 = X\cap P$, Conclusion (\ref{gopher}) of the
lemma follows.  Conclusions (\ref{reg nbhd}), (\ref{more reg nbhd})
and (\ref{twofer}) are immediate from the construction of $\calw$.
\EndProof

\Lemma\label{oinksday}
Suppose that $\calw$ is a  book
of $I$-bundles, and let $p$ denote the number of pages of
$\calw$. Then
 $$\rk H_2(|\calw|;\Z_2)\le p.$$
\EndLemma

\Proof
It is most natural to prove a very mild generalization: if
$\calw$ is a generalized book of $I$-bundles whose bindings
are all solid tori, and if $p$ denotes the number of pages of $\calw$,
then $\rk H_2(|\calw|;\Z_2)\le p$. We set $W=|\calw|$ and use
induction on $p$. If $p=0$ then the components of $W$ are solid tori
and hence $\rk H_2(W)=0$.  If $p>0$, choose a page $P$ of $\calw$, and
set $W'=\overline{W-P}$ and $ \calp'=\calp_\calw-P$. Then $\calp'$
inherits an $I$-bundle structure from $\calp$, and
$\calw'=(W',\calb,\calp')$ is a book of $I$-bundles with $p-1$ pages.
By the induction hypothesis we have $\rk H_2(W')\le p-1$. On the other
hand, if $F$ denotes the base surface of the $I$-bundle $P$, we have
$$H_2(W,W';\Z_2)\cong H_2(P,\partial_vP;\Z_2)\cong H_2(F,\partial
F;\Z_2)$$
and hence $\rk H_2(W,W')=1$. It follows that
$$\rk H_2(W)\le\rk H_2(W')+\rk H_2(W,W')\le p.$$
\EndProof

\Lemma\label{moosday} If $\calw$ is a  book of
$I$-bundles, and if every page of $\calw$ has strictly negative
Euler characteristic, we have
 $$\rk(|\calw|)\le2\barchi(|\calw|)+1.$$
\EndLemma

\Proof Set $W=|\calw|$. By hypothesis we have
$\chibar(P)\ge1$ for every page $P$ of $\calw$. Hence if
$P_1,\ldots,P_p$ denote the pages of $\calw$, we have
$$\chibar(W)=\sum_{i=1}^p\chibar(P_i)\ge p.$$
According to Lemma \ref{oinksday} we have
$$\rk H_2(W;\Z_2)\le p\le\chibar(W).$$
Now $W$ is a connected $3$-manifold with non-empty boundary. Hence
$\rk H_0(W;\Z_2)=1$, and $H_j(W;\Z_2)=0$ for each $j>2$. In view of \ref{oiler}, we have
$$\chibar(W)=\rk H_1(W;\Z_2)-\rk H_2(W;\Z_2)-1.$$
Hence
$$\rk(W)=\rk H_1(W;\Z_2)=\chibar(W)+\rk H_2(W;\Z_2)+1\le2\chibar(W)+1.$$
\EndProof

\section {Compressing submanifolds}
\label{compressing section}

\Definition\label{complexity def}
If $\cals$ is a closed (possibly empty or disconnected) surface, we define
a non-negative integer $\kappa(\cals)$ by
$$\kappa(V)=\sum_S(1+\genus(S)^2),$$
where $S$ ranges over the components of $\cals$.
\EndDefinition

\Lemma\label{pere tranquille}
Let $\cals$ be a closed (possibly empty or disconnected) surface, let
$A\subset \cals$ be a homotopically non-trivial annulus, and let
$\cals'$ be the surface obtained from the bounded surface
$\overline{\cals-A}$  by attaching disks $D_1$ and $D_2$ to
its two boundary components. Then $\kappa(\cals')<\kappa(\cals)$.
\EndLemma

\Proof
Let us index the components of $\cals$ as $S_0,\ldots,S_n$, where
$n\ge0$ and $A\subset S_0$. If $S_0-A$ is connected, the components of
$\cals'$ are $S_0',S_1,\ldots,S_n$, where $S_0'=(S_0-A)\cup
D_1\cup D_2$. We then have $\genus S_0'=(\genus S_0)-1$, so that $\kappa(
\cals)<\kappa(\cals')$. 

If $S_0-A$ is disconnected, then $(S_0-A)\cup D_1\cup D_2$ has two
components $S_0'$ and $S_0''$. If we denote the respective genera of
$S_0$, $S_0'$ and $S_0''$ by $g$, $g'$ and $g''$, we have $g=g'+g''$;
and since $A$ is homotopically non-trivial in $S_0$, both $g'$ and
$g''$ are strictly positive.  It follows that
$(1+(g')^2)+(1+(g'')^2)<1+g^2$, and we again deduce that $\kappa(
\cals)<\kappa(\cals')$.
\EndProof

\Number\label{des moines} Recall that a connected $3$-manifold $H$ is
called a {\it compression body} if it can be constructed from a
product $T\times[-1,1]$, where $T$ is a connected, closed, orientable
$2$-manifold, by attaching finitely many $2$- and $3$-handles to
$T\times\{-1\}$.  One defines $\partial_+H$ to be the submanifold
$T\times\{1\}$ of $\partial H$, and one define $\partial_-H$ to be
$\partial H-\partial_+ H$.  \EndNumber

\Number\label{waterloo}
If $H$ is a connected compression body, it is clear that $\partial_+H$
is connected, and that for each component $F$ of $\partial_-H$ we have
$\genus(F)\le\genus(\partial_+H)$. 
\EndNumber

\Number\label{cedar rapids}
It is a standard observation that a connected compression body $H$
with $\partial_-H=\emptyset$ is a handlebody.
\EndNumber

\Number\label{davenport}
Another standard observation is that any connected compression body
$H$ with $\partial_H\ne\emptyset$ can be constructed from a product
$S\times[-1,1]$, where $S$ is a possibly disconnected, closed,
orientable $2$-manifold, by attaching $1$- and $2$-handles to
$S\times\{1\}$. One then has $\partial_-H=S\times\{-1\}$.

An immediate consequence of this observation is that if $H$ is a
connected compression body then $\partial_-H$ is
$\pi_1$-injective in $H$.
\EndNumber

\Number\label{sioux city}
More generally, we shall define a {\it compression body} to be a
compact, possibly disconnected $3$-manifold $\calh$ such that each
component of $\calh$ is a compression body in the sense defined
above. We define  $\partial_+\calh=\bigcup_H\partial_+H$ and
$\partial_-\calh=\bigcup_H\partial_-H$, where $H$ ranges over the
components of $\calh$.
\EndNumber

\Proposition\label{new old prop 2}
Let $N$ be a compact orientable, irreducible $3$-manifold, and let $V$
be a compact, connected, non-empty $3$-submanifold of $\inter N$.
Suppose that $\overline{N-V}$ is $\pi_1$-injective in $N$. Then at
least one of the following conditions holds:
\Conclusions
\item\label{gone the way of brother tom}$V$ is contained in a ball in
  $N$; or 
\item\label{sister jenny's turn} $\partial V\ne\emptyset$, and there
  exists a connected, incompressible closed surface  in $N$ whose
  genus is at most the maximum of the genera of the
  components of $\partial V$; or
\item\label{mama's aim is bad} $N$ is closed and  every component of
  $\overline{N-V}$  is  a handlebody. 
\EndConclusions
\EndProposition

\Proof 
First note that if $V=N$ then conclusion (\ref{mama's aim is bad})
holds. (The hypothesis $V\subset\inter N$ implies that $N$ is closed,
and the other assertion of (\ref{mama's aim is bad}) is vacuously
true.) Hence we may assume that $V\ne N$, so that $\partial
V\ne\emptyset$.

Let $\calc$ denote the set of all (possibly disconnected) compression
bodies $\calh\subset N$ such that $\calh\cap
V=\partial_+\calh=\partial V$. Note that a regular neighborhood of
$\partial V$ relative to $\overline{N-V}$ is an element of $\calc$,
and hence that $\calc\ne\emptyset$. Let us fix an element $\calh$ of
$\calc$ such that (in the notation of \ref{complexity def}) we have
$\kappa(\partial_-\calh)\le\kappa(\partial_-\calh')$ for every
$\calh'\in\calc$.

Note that $V\cup \calh$ is connected since $\calh\in\calc$.

Consider first the case in which $\partial_- \calh=\emptyset$. In this
case, it follows from \ref{cedar rapids} that every component of $\calh$
is a handlebody, and we have $\partial \calh=\partial_+\calh=\partial
V$. Since $N$ is connected and $\partial V\ne\emptyset$, we must have
$\calh=\overline{N-V}$. In particular $N$ must be closed. Thus conclusion
(\ref{mama's aim is bad}) of the proposition holds in this case.

Now consider the case in which some component $S$ of $\partial_-\calh$
is a $2$-sphere. By irreducibility, $S$ bounds a ball $B\subset
N$. Since $V\cup \calh$ is connected, we have either $V\cup
\calh\subset B$ or $B\cap(V\cup \calh)=\partial B$. If $V\cup
\calh\subset B$, then in particular conclusion (\ref{gone the way of
  brother tom}) of the proposition holds. If $B\cap(V\cup
\calh)=\partial B$, then $\calh'\doteq \calh\cup B$ is obtained from
$\calh$ by attaching a $3$-handle to $\partial_-\calh$, and hence
$\calh'\in\calc$ (cf. \ref{des moines}). But we have
$\partial_-\calh'=\partial_-\calh-S$, and it follows from Definition
\ref{complexity def} that $\kappa(\calh')=\kappa(\calh)-1$. This
contradicts the minimality of $\kappa(\calh)$.

There remains the case in which $\partial_- \calh\ne\emptyset$, and every
component of $\partial_- \calh$ has positive genus. Let us fix a component
$Z$ of $\overline{N-V}$ which contains at least one component of
$\partial_-\calh$. Let us set $F=Z\cap\partial_-\calh$. Then $F$ is a
non-empty (and possibly disconnected) closed surface in $\inter Z$,
and each component of $F$ has positive genus. We claim that $F$ is
incompressible in $Z$.

Suppose to the contrary that $F$ is compressible in $Z$. Then there is
a disk $D\subset\inter Z$ such that $D\cap F=\partial D$, and such
that $\partial D$ is a homotopically non-trivial simple closed curve
in $F$.  Since $D\subset\inter Z\subset N-V$, we have
$D\cap\partial_+\calh=\emptyset$. Furthermore, since $D\subset Z$, we have
$D\cap\partial_-\calh=D\cap( Z\cap\partial_-\calh)=D\cap F=\partial D$. Hence
$D\cap\partial \calh=\partial D$. It follows that either $D\subset \calh$ or
$D\cap \calh=\partial D$. 

If $D\subset \calh$, let $H_0$ denote the component of $\calh$ containing $D$,
and let $F_0\subset\partial_-H_0$ denote the component of $F$
containing $\partial D$. Since $\partial D$ is homotopically
non-trivial it follows that the inclusion homomorphism
$\pi_1(F_0)\to\pi_1(H_0)$ has non-trivial kernel. This contradicts
\ref{davenport}. 

If $D\cap \calh=\partial D$, we fix a regular neighborhood $E$ of $D$
relative to $\overline{Z-\calh}$. Then $\calh'\doteq \calh\cup E$ is
obtained from $\calh$ by attaching a $2$-handle to $\partial_-\calh$,
and hence $\calh'\in\calc$ (cf. \ref{des moines}). The surface
$\partial \calh'$ has the form $((\partial \calh)-A)\cup D_1\cup D_2$,
where $A\subset\partial \calh$ is a homotopically non-trivial annulus,
and $D_1$ and $D_2$ are disjoint disks in $N$ such that $(D_1\cup
D_2)\cap\partial \calh=\partial A$. It therefore follows from Lemma
\ref{pere tranquille} that $\kappa( \partial \calh)<\kappa(\partial
\calh')$. This contradicts the minimality of $\kappa(\calh)$, and the
incompressibility of $F$ in $Z$ is proved.

Since $Z$ is $\pi_1$-injective in $N$ by hypothesis, it now follows
that $F$ is incompressible in $N$. Our choice of $Z$ guarantees that
$F\ne\emptyset$. Choose any component $F_1$ of $F$, and let $H_1$
denote the component of $\calh$ containing $F_1$. By \ref{waterloo},
$\genus(F_1)$ is at most the genus of the connected surface
$\partial_+ \calh$. But $\partial_+\calh$ is a component of $\partial
V$ since $\calh\in\calc$, and so $\genus(F_1)$ is at most the maximum
of the genera of the components of $\partial V$. Hence conclusion
(\ref{sister jenny's turn}) of the proposition holds in this case.
\EndProof

\section{Transporting surfaces downstairs}
\label{down}

\Lemma\label{the chalice in the palace}
Let $M$ be a simple, compact, orientable $3$-manifold, let $p:\tM\to
M$ be a $2$-sheeted covering, and let $\tau:\tM\to\tM$ denote the
non-trivial deck transformation.  Suppose that $\tM$ contains a
closed, incompressible surface $F_0$ of positive genus.  Then $F_0$ is
ambiently isotopic to a surface $F$ such that $F$ and $\tau(F)$ meet
transversally, and every component of $F\cap\tau(F)$ is a
homotopically non-trivial simple closed curve in $M$.
\EndLemma

\Proof
Let $\calf$ denote the collection of all surfaces $S\subset M$ such that
$S$ is isotopic to $F_0$ and $S$ meets $\tau(S)$ transversely. 
Choose a surface $F\in\calf$ so that the number of components of
$F\cap\tau(F)$ is minimal.  We will show that every component
of $F\cap\tau(F)$ is a homotopically non-trivial simple closed curve.

Suppose there exists a homotopically trivial component $\gamma$ of
$F\cap\tau(F)$.  Then, since $F$ is incompressible in $M$, the simple
closed curve $\gamma$ must bound disks $D\subset F$ and $D'\subset
\tau(F)$.  We assume, without loss of generality, that the disk $D'$
is innermost on $\tau(F)$ in the sense that $ D'\cap F = \gamma$.
This implies, in particular that $D\cup D'$ is an embedded $2$-sphere
in $M$.

 Since $M$ is irreducible, the $2$-sphere $D\cup D'$ bounds a ball $B$
in $M$.  We may observe at this point that the curve $\gamma$ cannot
be invariant under $\tau$.  Otherwise, since $D'$ is the unique disk
on $\tau(F)$ bounded by $\gamma$, it would follow that $\tau(D)=D'$,
and hence that the sphere $D\cup D'$ is invariant under $\tau$.  Since
$\tM$ contains an incompressible surface, it is not homeomorphic to
$S^3$, and therefore $B$ is the unique $3$-ball bounded by $D\cup D'$.
Thus the assumption that $\gamma$ is invariant implies that the ball
$B$ is invariant under the fixed point free map $\tau$, contradicting
the Brouwer Fixed Point Theorem.  This shows that $\gamma$ is not
invariant under $\tau$.  It follows, since $D'$ is innermost, that
$D'$ is disjoint from its image under $\tau$.

Now let $V$ be a regular neighborhood of $B$, chosen so that $V\cap F$
is a regular neighborhood of $D$ and $V\cap F'$ is a regular
neighborhood of $D'$.  The disk $F'\cap V$ divides $V$ into two balls,
one of which, say $U$, is disjoint from the interior of $D$.  Since
$D'\cap \tau(D') = \emptyset$, we may assume without loss of
generality that $V$ has been chosen to be small enough so that
$U\cap\tau(U) = \emptyset$.  Let $E$ denote the disk in $\partial U$
which is bounded by $F\cap U$ and which is disjoint from $\tau(F)$.
We set $A =\overline{F\setminus U}$ and consider the surface $F' =
A\cup E$, which is clearly isotopic to $F$ by an isotopy supported in
$V$.  We will show that $F'\cap\tau(F') \subset (F\cap \tau(F)) -
\gamma$.

We write $F'\cap\tau(F') = (A\cup E)\cap(\tau(A)\cup\tau(F))$ as the
union of the four sets $A\cap\tau(A)$, $A\cap\tau(E)$, $E\cap\tau(A)$
and $E\cap\tau(E)$.  We have $A\cap\tau(A)\subset
F\cap\tau(F)-\gamma$.  Since $E\subset U$ and
$U\cap\tau(U)=\emptyset$ we have $E\cap\tau(E)=\emptyset$.  The
sets $E$ and $\tau(F)\supset\tau(A)$ are disjoint by construction, and
hence $E\cap\tau(A) = \emptyset$. Finally,
$A\cap\tau(E)=\tau(E\cap\tau(A)) = \emptyset$.

We have shown that $F'\cap\tau(F') \subset (F\cap \tau(F)) - \gamma$,
and hence that  $F'\cap\tau(F')$ has fewer components that
$F\cap \tau(F)$.  This contradicts the choice of $F$, and completes
the proof of the lemma.
\EndProof

\Lemma\label{the brew that is true} Let $ N $ be a simple, compact,
orientable $3$-manifold, let $p: \tN \to N $ be a $2$-sheeted
covering, and let $\tau: \tN \to \tN $ denote the non-trivial deck
transformation.  Suppose that $F\subset \tN $ is a closed,
incompressible surface such that $F$ and $\tau(F)$ meet transversally,
and every component of $F\cap\tau(F)$ is a homotopically non-trivial
simple closed curve in $ N $. Then $N -p(F)$ is $\pi_1$-injective
in $N$.
\EndLemma

\Proof
Set $F_1=\tau(F)$, so that $F_1$ is incompressible in $\tN$. Set
$C=F\cap F_1$.  Let $\tN'$ denote the $3$-manifold obtained by splitting $\tN$
along $F$, and let $F_1'$ denote the surface obtained by splitting
$F_1$ along $C$. Then $\tN$ and $F_1$ may be regarded as quotient
spaces of $\tN'$ and $F_1'$, and 
$F_1'$ is naturally identified with a properly embedded surface in $\tN'$.
We have a commutative diagram
$$\xymatrix{
F_1'   \ar[r] \ar[d]       &   F_1   \ar[d] \\
\tN'   \ar[r]              &   \tN \\
}$$
where the horizontal maps are quotient maps and the vertical maps are
inclusions. The inclusion $F_1\to\tN$ is $\pi_1$-injective because
$F_1$ is incompressible in $\tN$, and the quotient map $F_1'\to F_1$
is $\pi_1$-injective because the components of $C$ are homotopically
non-trivial. By commutativity of the diagram it follows that the
inclusion $F_1'\to\tN'$ is $\pi_1$-injective.

Now let $\tN''$ denote the $3$-manifold obtained by splitting $\tN'$
along $F_1'$. Since the
inclusion $F_1'\to\tN'$ is $\pi_1$-injective, the quotient map
$\tN''\to\tN'$ is also $\pi_1$-injective. On the other hand, the quotient map
$\tN'\to\tN$ is $\pi_1$-injective because $F$ is incompressible in
$\tN$. Hence the composite quotient map $\tN''\to\tN$ is
$\pi_1$-injective. It follows that the inclusion map $\tN-(F\cup
F_1)\to\tN$ is $\pi_1$-injective.

Now consider any component $Z$ of $ N -p(F)$.  Choose a component
$\tZ$ of $p^{-1}(Z)$. Then $\tZ$ is a component of $\tN-(F\cup F_1)$,
and hence the inclusion $\tZ\to\tN$ is $\pi_1$-injective.  Thus in the
commutative diagram
$$\xymatrix{
\pi_1(\tZ) \ar[r] \ar[d]   &   \pi_1(\tN) \ar[d] \\
\pi_1(Z)   \ar[r]          &   \pi_1(N) \\
}$$
the inclusion homomorphism $\pi_1(\tZ)\to \pi_1(\tN)$ is injective,
while the vertical homomorphisms are induced by covering maps and are
therefore also injective. Since the image of $\pi_1(\tZ)$ has index at
most $2$ in $\pi_1(Z)$, the kernel of the inclusion homomorphism
$\pi_1(Z)\to \pi_1(N)$ has order at most $2$. But $\pi_1(Z)$ is
torsion-free because $N$ is simple.  Hence
$\pi_1(Z)\to \pi_1(N)$ is injective, as asserted by the Lemma.
\EndProof

\Lemma\label{new lemma for an old prop} Suppose that $  N  $ is a
simple, compact, orientable $3$-manifold, that $p:  \tN  \to   N  $ is a
$2$-sheeted covering, that $g\ge2$ is an integer, and that $  \tN  $
contains a closed, incompressible surface of genus $g$.  Then there exist a
 connected
   book of $I$-bundles $\calv$ with $V=|\calv|\subset   N  $, and a 
 closed, orientable (possibly disconnected) surface $S\subset \inter
 V$ such that 
\Conclusions
\item\label{you can go by foot}$\chibar(V)=\chibar(S)=2g-2$;
\item\label{you can go by cow}every page of $\calv$ has strictly negative
 Euler characteristic;
\item\label{i'm injective}$\calp_\calw$ is $\pi_1$-injective in $N$;
\item\label{you're injective}$ N -V$ is
$\pi_1$-injective in $N$; 
\item\label{square hole}no component of $S$ is a sphere; and
\item\label{marvin k mooney}for every page $P$ of
 $\calv$, the set $S\cap P$ is a section of the $I$-bundle $P$.
\EndConclusions
\EndLemma

\Proof 
According to Lemma \ref{the chalice in the palace}, $ \tN $ contains a
closed, incompressible surface $F$ of genus $g$ such that $F$ and
$\tau(F)$ meet transversally, and every component of $F\cap\tau(F)$ is
a homotopically non-trivial simple closed curve in $ N $.  It
follows that $q=p|F:F\to N$ is an immersion with at most double-curve
singularities. 

The map $q_\sharp:\pi_1(F)\to\pi_1(N)$ is injective because $F$ is
incompressible in $\tN$ and $p:\tN\to N$ is a covering map.

Let us set $X=q(F)$, and let $C\subset X$ denote the union of all
double curves of $q$.  Since the components of $C$ are homotopically
non-trivial in $N$ and hence in $X$, the set $\tC=q^{-1}(C)$ is a
disjoint union of homotopically non-trivial simple closed curves in
$F$. Hence $F-\tC$ is $\pi_1$-injective in $F$, and each of its
components has non-positive Euler characteristic.  Since
$q_\sharp:\pi_1(F)\to\pi_1(N)$ is injective it follows that
$q|(F-\tC):(F-\tC)\to N$ is $\pi_1$-injective.

The set $F-\tC$ is mapped homeomorphically onto $X-C$ by $q$. In
particular, each component of $X-C$ has non-positive Euler
characteristic. Furthermore, since $q|(F-\tC):(F-\tC)\to N$ is
$\pi_1$-injective, it now follows that $X-C$ is $\pi_1$-injective in
$N$.

In the notation of (\ref{Pi notation}), we have $\Pi(X)=X-C$ , and the
link in $X$ of every point of $C$ is homeomorphic to the suspension of
a four-point set.  Since every component of $X-C$ has non-positive
Euler characteristic, it follows from Definition \ref{bosdef} that $X$
is a book of surfaces.  Since each component of $C$ is a simple closed
curve, we have $\barchi(X)=\barchi(F)=2g-2$.

Let $W$ denote a regular neighborhood of $X$ in $ N $. According to
Lemma \ref{neighborhood is a book}, we may write $W=|\calw|$ for some
book of $I$-bundles $\calw$ in such a way that  Conclusions
(\ref{more reg nbhd})--(\ref{twofer}) of Lemma \ref{neighborhood is a
book} hold. Since $X-C$
is $\pi_1$-injective in $N$, it follows from Conclusion
(\ref{twofer}) of Lemma \ref{neighborhood is a book} that
$\calp_\calw$ is $\pi_1$-injective in $N$.

Since $\chi(W)=\chi(X)=2-2g<0$, and since $\calp_\calw$ is
$\pi_1$-injective in $N$, it follows from Lemma \ref{make book}
that there is a connected book of $I$-bundles $\calv$ with
$V=|\calv|\subset N $, such that Conclusions (\ref{a pair of pizza
pies})---(\ref{mutilated monkey meat}) of Lemma \ref{make book} hold.
Conclusion (\ref{Euler doesn't care}) of Lemma \ref{make book} gives
$\chibar(V)=\chibar(W)=\chibar(X)$, so that \Equation\label{a high
wind in the attic} \chibar(V)=2g-2.  \EndEquation 

It follows from Conclusions (\ref{a pair of pizza pies}) and (\ref{and
me without a spoon}) of Lemma \ref{make book} that every binding of
$\calw$ is contained in a binding of $\calv$. Since by Conclusion
(\ref{more reg nbhd}) of Lemma \ref{neighborhood is a book} we have
$C\subset\inter\calb_\calw$, it follows that
$C\subset\inter\calb_\calv$.

Let $\calu$ denote a regular neighborhood of $C$ in $\inter\calb_\calv$. We may
suppose $\calu$ to be chosen so that $\partial \calu$ meets $\Pi(X)$
transversally, and each component of $\calu\cap X$ is homeomorphic to
$\text{\large +}\times S^1$, where $\text{\large +}$ denotes a cone on
a four-point set. Set $X'=\overline{X-(\calu\cap X)}$ and $F'=F\cap
q^{-1}(X')$. Then $F'$ and $X'$ are (possibly disconnected) compact
$2$-manifolds with boundary, and $q'=q|F'$ maps $F'$ homeomorphically
onto $X'$. Let us fix an orientation of $F$, so that $F'$ inherits an
orientation, and define an orientation of $X'$ by transporting the
orientation of $F'$ via $q$.

Let $U_1,\ldots,U_m$ denote the components of $\calu$.  We set
$B_i=X\cap\partial U_i$. Each component $\beta$ of $B_i$ is a boundary
component of $X'$ and hence has an orientation induced from the
orientation of $X'$, which determines a generator of $H_1(U_i;\Z)$ via
the inclusion isomorphism $H_1(\beta;\Z)\to H_1(U_i;\Z)$. We shall say
that two components of $B_i$ are {\it similar} if they determine the
same generator of $H_1(U_i;\Z)$ via this construction.

The set $(\partial U_i)-B_i$ has four components. Their closures are
annuli, which we shall call {\it complementary annuli}.  We shall say
that two components of $ B_i$ are {\it adjacent} if their union
is the boundary of a complementary annulus, and {\it opposite}
otherwise.

If $\beta$ and $\beta'$ are opposite components of $ X\cap\partial
U_i$, then $q^{-1}(\beta)$ and $q^{-1}(\beta')$ form the boundary
of an annulus $A$ in $F$, which is mapped homeomorphically by
$q$ to an embedded annulus in $U_i$. Since the orientation of $F'$
is the restriction of an orientation of $F$, the induced
orientations of $q^{-1}(\beta)$ and $q^{-1}(\beta')$ determine
different generators of $H_1(A;\Z)$. In view of our definitions it
follows that opposite components of $ B_i$ are dissimilar.

Let us call a complementary annulus {\it bad} if its boundary
curves are similar, and {\it good} otherwise.  If $\beta$ is any
component of $B_i$, the two components of $B_i$ adjacent to $\beta$
are opposite each other; hence exactly one of them is similar to
$\beta$. This shows that $\beta$ is contained in the boundary of
exactly one bad annulus and one good annulus. We conclude that
$\partial U_i$ contains exactly two good annuli, say $A_i$ and $A_i'$,
and that $A_i\cap A_i'=\emptyset$.

The set
$$S=(X-(X\cap\calu))\cup(A_1\cup\cdots\cup A_m)\cup(A'_1\cup\cdots\cup
A'_m)$$ is a (possibly disconnected) compact PL
$2$-manifold embedded in $V$. Since $A_i$ and $A_i'$ are good annuli,
the orientation of $X'$ extends to an orientation of $S$. In particular
$S$ is orientable. 

We shall show that Conclusions (\ref{you can go by foot})--(\ref
{marvin k mooney}) of the present lemma hold when $\calv$
and $S$ are defined as above.

According to (\ref{a high
wind in the attic}) we have $\chibar(V)=2g-2$. It follows from the
construction of $S$ that $\chibar(S)=\chibar(X)=2g-2$. Hence Conclusion
(\ref{you can go by foot}) of the present lemma holds. 

Conclusion (\ref{you can go by cow}) of the present lemma follows from
 Conclusion (\ref{he does and he doesn't}) of Lemma \ref{make book}.

Since we have seen that $\calp_\calw$ is $\pi_1$-injective in $N$, it
follows from Conclusion (\ref{and me without a spoon}) of Lemma
\ref{make book} that $\calp_\calv$ is $\pi_1$-injective in $N$. This
is Conclusion (\ref{i'm injective}) of the present lemma.

 It follows from Lemma \ref{the brew that is true} that $N-X= N-q(F)$
is $\pi_1$-injective in $N$. It follows from conclusions (\ref{a pair
of pizza pies}) and (\ref{wysiwyg}) of Lemma \ref{make book} that every
component of $N-V$ is also a component of $N-W$, and is therefore
ambiently isotopic in $N$ to a component of $ N- X$. Hence $ N -V$ is
$\pi_1$-injective in $N$. This is Conclusion (\ref{you're injective})
of the present lemma.

It follows from the construction of $S$ that
$S\cap\calp_\calv=X\cap\calp_\calv$. If $P$ is any page of $\calv$,
then by Conclusion (\ref{and me without a spoon}) of Lemma \ref{make
  book}, $P$ is a page of $\calw$, and hence $S\cap P=X\cap P$ is a
section of the $I$-bundle $P$ according to Conclusion (\ref{gopher})
of Lemma \ref{neighborhood is a book}. This establishes Conclusion
(\ref{marvin k mooney}) of the present lemma.

In particular it follows that for every page $P$ of $\calv$ the
surface $P\cap S$ is connected and has non-positive Euler
characteristic. On the other hand, the construction of $S$ shows that
every component of $S\cap\calb_\calv$ is an annulus. Hence every
component of $S$ has non-positive Euler characteristic, and Conclusion
(\ref{square hole}) of the present lemma follows.
\EndProof

\Proposition\label{new improved old prop 3} Suppose that $  N  $ is a
simple, compact, orientable $3$-manifold, that $p:  \tN  \to   N  $ is a
$2$-sheeted covering, that $g\ge2$ is an integer, and that $  \tN  $
contains a closed, incompressible surface of genus $g$.  Then either
\Conclusions
\item\label{hit 'im wid a brick} $N$ contains a closed, connected,
  incompressible surface of genus at most $g$, or
\item\label{bust 'em clown} $N$ is closed and there is a connected,
 book of $I$-bundles $\calv$ with $V=|\calv|\subset N $ such that
 $\chibar(V)=2g-2$, every page of $\calv$ has strictly negative
 Euler characteristic, and every component of $\overline{N-V}$ is a
 handlebody.  In particular, the rank of $H_1(N;\Z_2)$ is at most
 $4g-3$.  Furthermore, there is a closed, orientable (possibly
 disconnected) surface
 $S\subset \inter V$ such that for every
 page $P$ of $\calv$, the set $S\cap P$ is a section of the $I$-bundle
 $P$.
\EndConclusions 
\EndProposition

The last sentence of alternative of (\ref{bust 'em clown}) is not
used in this paper, but will be needed in \cite{second}.

\Proof[Proof of Proposition \ref{new improved old prop 3}]
Let us fix a connected book of $I$-bundles $\calv$ with
$V=|\calv|\subset N $, and a closed, orientable surface $S\subset
\inter V$, such that Conclusions
(\ref{you can go by foot}) to (\ref{marvin k mooney}) of Lemma \ref{new
lemma for an old prop} hold.  We distinguish two cases, depending on
whether there does or does not exist a page of $\calv$ whose
horizontal boundary is contained in a single component of $\partial V$.

{\bf Case I. There is a page $P_0$ of $\calv$ such that $\partial_hP_0$ is
contained in a single component $Y_0$ of $\partial V$.}

According to conclusion (\ref{marvin k mooney}) of Lemma \ref{new
  lemma for an old prop}, the set $S\cap P_0$ is a section of the
$I$-bundle $P_0$.  Hence there is a properly embedded arc $\alpha$ in
$V$, such that $\alpha\subset P_0$, and such that $\alpha$ meets $S$
transversally in a single point. The endpoints of $\alpha$ lie in
$\partial_hP_0\subset Y_0$. Since $Y_0$ is connected, there is an arc
$\beta\subset Y_0$ with $\partial\beta=\partial\alpha$.

Let $\sigma$ denote the class in $ H_2(N;\Z_2)$ represented by
$S$. Since $\alpha$ is properly embedded in $V$ and meets $X$
transversally in a single point of $\pi_0\subset\Pi(X)$, the class
$\sigma$ has intersection number $1$ with the class in $H_1(N;\Z_2)$
represented by the simple closed curve $\alpha\cup\beta$. In
particular $\sigma\ne0$.  Hence some component $S_0$ of $S$ represents
a non-zero class in $H_2(N;\Z_2)$. It follows from Conclusions
(\ref{you can go by foot}) and (\ref{square hole}) of Lemma \ref{new
  lemma for an old prop} that $\chibar(S_0)\le\chibar(S)=2g-2$, and
hence that $\genus(S_0)\le g$.

Among all closed, orientable surfaces in $ N $ that represent
non-trivial classes in $ H_2( N ;\Z_2)$, let us choose one, say $S_1$,
of minimal genus. Then $\genus (S_1)\le\genus(S_0)\le g$.  If $S_1$ is
compressible in $ N $, a compression of $S_1$ produces a $2$-manifold
$S_2$ with one or two components. Each component of $S_2$ has strictly
smaller genus than $S_1$, and at least one of them represents a
non-trivial class in $ H_2( N ;\Z_2)$.  This contradicts
minimality. Hence $S_1$ is incompressible in $ N $. Since
$\genus(S_1)\le g$, conclusion (\ref{hit 'im wid a brick}) of
the present lemma holds in this case.

{\bf Case II.  There is no page $P_0$ of $\calv$ such that $\partial_hP_0$ is
contained in a single component of $\partial V$.}

In this case, every page of $\calv$ is a trivial
$I$-bundle. Furthermore, if  $T$ is any component  of $\partial V$,
then for every page $P$ of $\calv$, at most one component of the
horizontal boundary of $P$ is contained in $T$. Hence 
$$\chibar(T\cap P)\le\chibar (P)$$ for every page $P$ of
$\calv$. Letting $P$ range over the pages of $\calv$, and using
(\ref{a high wind in the attic}), we find that
$$\chibar(T)=\sum_P\chibar(T\cap P)\le\sum_P\chibar(P)=\chibar(V)=2g-2.$$
This shows that
\Equation\label{rats in jamaica}
\genus(T)\le g
\EndEquation
for every component $T$ of $\partial V$.

According to Conclusion (\ref{you're injective}) of Lemma \ref{new
  lemma for an old prop}, $ N -V$ is
$\pi_1$-injective in $N$. Thus $V\subset N$ satisfies the hypotheses
of Proposition \ref{new old prop 2}. There are three subcases
corresponding to the three alternatives (\ref{gone the way of brother
tom})---(\ref{mama's aim is bad}) of Proposition \ref{new old prop 2}.

First suppose that alternative (\ref{gone the way of brother tom}) of
Proposition \ref{new old prop 2} holds, i.e. that $V$ is contained
in a ball. Then in particular for any page $P$ of $\calv$, the
inclusion homomorphism $\pi_1(P)\to\pi_1(W)$ is trivial.  But
according to Conclusions (\ref{you can go by cow}) and (\ref{i'm
  injective}) of Lemma \ref{new lemma for an old prop}, we have
$\chi(P)<0$ (so that $\pi_1(P)$ is non-trivial) and $\calp_\calw$ is
$\pi_1$-injective in $N$. This contradiction shows that alternative
(\ref{gone the way of brother tom}) of Proposition \ref{new old prop
  2} cannot hold in this situation.

Next suppose that alternative (\ref{sister jenny's turn}) of
Proposition \ref{new old prop 2} holds, i.e. that there exists a
connected, incompressible closed surface $S_1$ in $ N $ whose genus is
at most the maximum of the genera of the components of $\partial
V$.  By (\ref{rats in jamaica}) this maximum is at most $g$.
Thus conclusion (\ref{hit 'im wid a brick}) of the present lemma holds
in this subcase.

Finally, suppose that alternative (\ref{mama's aim is bad}) of
Proposition \ref{new old prop 2} holds, i.e. that $N$ is closed and
that every component of $\overline{N-V}$ is a handlebody.  We have
that $V=|\calv|$ where $V$ is a book of $I$-bundles whose pages
all have negative Euler characteristic, and Conclusion (\ref{you can
go by foot}) of Lemma \ref{new lemma for an old prop} gives
$\chibar(V)= 2g-2$. Since the components of $\overline{N-V}$ are
handlebodies, the inclusion of $V$ into $N$ induces a surjection from
$H_1(V;\Z_2)$ to $H_1(N;\Z_2)$; hence the latter group has rank at
most $4g-3$ by Lemma \ref{moosday}.  Furthermore, according to
Conclusion (\ref{marvin k mooney}) of Lemma \ref{new lemma for an old
prop}, for every page $P$ of $\calv$, the set $S\cap P$ is a section
of the $I$-bundle $P$.  Thus conclusion (\ref{bust 'em clown}) of the
present proposition holds in this subcase.  \EndProof

\section{Singularity of PL maps}
\label{singularity section}

If $K$ is a finite simplicial complex, we shall denote the
underlying space of $K$ by $|K|$. A simplicial map $\phi:K_1\to K_2$ 
between finite simplicial complexes defines a map from $|K_1|$ to $|K_2|$
which we shall denote by $|\phi|$.

Now suppose that $X_1$ and $X_2$ are compact topological spaces and
that $f:X_1\to X_2$ is a continuous surjection. We define a {\it
  triangulation} of $f$ to be a quintuple $(K_1,J_1,K_2,J_2,\phi)$,
where each $K_i$ is a finite simpicial complex, $J_i:|K_i|\to X_i$ is
a homeomorphism, and $f\circ J_1=J_2\circ\phi$. When it is unnecessary
to specify the $K_i$ and $J_i$ we shall simply say that $\phi$ is a
triangulation of $f$.

Note that if $f$ is any PL map from a compact PL space $X$ to a PL
space $Y$, then the surjection $f:X\to f(X)$ admits a triangulation.

\Definition 
Let $K$ and $L$ be finite simplicial complexes and let $\phi:K\to L$
be a simplicial map. We define the {\it degree of singularity} of
$\phi$, denoted $\DS(\phi)$, to be the number of ordered pairs $(v,w)$
of vertices of $K$ such that $v\ne w$ but $\phi(v)=\phi(w)$.

If $f$ is any PL map from a compact PL space $X$ to a PL space $Y$,
we define the {\it absolute degree of singularity} of $f$, denoted
$\ADS(f)$, by
$$\ADS(f)=\min_{\phi}\DS(\phi),$$
where $\phi$ ranges over all
triangulations of $f:X\to f(X)$.
\EndDefinition

\Number\label{co-restriction}
We emphasize that the definition of $\ADS(f)$ is based on regarding
$f$ as a map from $X$ to $f(X)$. Hence if  $f$ is any PL map from a
compact PL space $X$ to a PL space $Y$, and $Z$ is a PL subspace of
$Y$ containing $f(X)$, then the absolute degree of singularity of $f$
is unchanged when we regard $f$ as a PL map from $X$ to $Z$.
\EndNumber

An almost equally trivial immediate consequence of the definition of
absolute degree of singularity is expressed by the following result.

\Lemma\label{homeomorphic images}
Suppose that $X$, $Y$ and $Z$ are PL spaces, that $X$ is compact, that
$f:X\to Y$ is a PL map, and that $h$ is a PL homeomorphism of $f(X)$
onto a PL subspace of $Z$. Then $h\circ f:X\to Z$ has the same
absolute degree of singularity as $f$.
\EndLemma

\Proof
In view of \ref{co-restriction} we may assume that $f$ is surjective
and that $h$ is a PL homeomorphism of $Y$ onto $Z$. Now if
$(K_1,J_1,K_2,J_2,\phi)$ is a triangulation of $f$ then $(K_1,J_1,
K_2,h\circ J_2,h\circ\phi)$ is a triangulation of $h\circ f$, and
$\DS(h\circ\phi)=\DS(\phi)$. It follows that $\ADS(h\circ
f)\le\ADS(f)$. The same argument, with $h^{-1}$ in place of $h$, shows
that $\ADS( f)\le\ADS( h\circ f)$.
\EndProof

\Proposition[Stallings]\label{stallings}
Suppose that $Y$ is a connected PL space and that $p:\widetilde Y\to
Y$ is a connected covering space, which is non-trivial in the sense
that $p$ is not a homeomorphism. Suppose that $f$ is a PL map from a
compact connected PL space $X$ to $Y$, such that the inclusion
homomorphism $\pi_1(f(X))\to\pi_1(Y)$ is surjective. Suppose that
$\widetilde f:X\to \widetilde Y$ is a lift of $f$. Then
$\ADS(\widetilde f)<\ADS(f)$.
\EndProposition

\Proof
We first prove the proposition in the special case where $f:X\to Y$ is
a surjection. In this case we set $m=\ADS(f)$, and we fix a
triangulation $(K_1,J_1,K_2,J_2,\phi)$ of the PL surjection $f$ such
that $\DS(\phi)=m$. Here, by definition, $J_1:|K_1|\to X$ and
$J_2:|K_2|\to Y$ are homeomorphisms. Let us identify $X$ and $Z$ with
$|K_1|$ and $|K_2|$ via these homeomorphisms. The covering space
$\widetilde Y$ of $Y$ may be identified with $|\widetilde K_2|$ for
some simplicial covering complex $\widetilde K_2$ of $K_2$; thus
$p=|\sigma|$ for some simplicial covering map $\sigma:\widetilde
K_2\to K_2$. The lift $\widetilde f$ may be written as
$|\widetilde\phi|$ for some simplicial lift
$\widetilde\phi:K_1\to\widetilde K_2$. We shall denote by $W$ the
subcomplex {$\widetilde\phi(K_1)$} of $\widetilde K_2$.

Since $\sigma\circ\widetilde\phi=\phi$, the definition of degree of
singularity implies that $\DS(\widetilde\phi)\le\DS(\phi)=m$. If
equality holds here, then the restriction of $\sigma$ to the vertex
set of $W$ is one-to-one. This implies that $p$ restricts to a
one-to-one map from $|W|$ to $Y$. But we have $W=\widetilde f(X)$, and
the surjectivity of $f$ implies that $p$ maps $|W|$ onto $Y$; thus $p$
restricts to a homeomorphism from $|W|$ to $Y$. This is impossible
since $p:\widetilde Y\to Y$ is a non-trivial connected covering space.
Hence we must have $\DS(\widetilde\phi)<m$. Since by definition we
have $\ADS(\widetilde\phi)\le DS(\widetilde\phi)$, the assertion of
the proposition follows in the case where $f$ is surjective.

We now turn to the general case. Let us set $Z=f(X)$ and $\widetilde
Z=p^{-1}(Z)$. Since $\widetilde Y$ is a non-trivial connected covering
space of $Y$, and since the inclusion homomorphism
$\pi_1(Z)\to\pi_1(Y)$ is surjective, $\widetilde Z$ is a non-trivial
connected covering space of $Z$. According to \ref{co-restriction},
regarding $\widetilde f$ and $f$ as maps into $\widetilde Z$ and $Z$ does not
affect their absolute degrees of  singularity. Since $f:X\to Z$ is
surjective,  the inequality now follows from the special case that has
already been proved.
\EndProof

Following the terminology used by Simon in \cite{simon}, we shall say
that a $3$-manifold $M$ admits a {\it manifold compactification} if
there is a homeomorphism $h$ of $M$ onto an open subset of a compact
$3$-manifold $Q$ such that $h(\inter M)=\inter Q$.

\Lemma\label{Simon}
Suppose that $N$ is a compact, orientable, connected, irreducible PL
$3$-manifold and that $D$ is a separating, properly embedded disk in
$N$. Let $X$ denote the closure of one of the connected components of
$N-D$. Let $\nu\in D$ be a base point, and let $p:(\widetilde
N,\widetilde\nu)\to (N,\nu)$ denote the based covering space
corresponding to the subgroup $\image(\pi_1(X,\nu)\to\pi_1(N,\nu))$ of
$\pi_1(N,\nu)$. Then $\widetilde N$ admits a manifold
compactification.
\EndLemma

\Proof
Let us set $X_1=\overline{N-X}$. It will also be convenient to
write $X_0=X$. Then the $X_i$ are compact submanifolds of $N$, and
$X_1\cap X_2=D$. We set $H_i=\pi_1(X_i,\nu)$ for $i=0,1$. We identify
$\pi_1(N,\nu)$ with $H_0\star H_1$, so that the $H_i$ are in
particular subgroups of $\pi_1(N,\nu)$.  Thus $(\widetilde N,\widetilde\nu)$
is the based covering space corresponding to the subgroup $H_0$.

According to the general criterion given by Simon in \cite[Theorem
3.1]{simon}, $\widetilde N$ will admit a manifold compactification
provided that the following conditions hold: 
\Conditions
\item $X_0$ and $X_1$ are irreducible,
\item $D$ is $\pi_1$-injective in $X_0$ and $X_1$,
\item each conjugate of $H_0$ in $\pi_1(N,\nu)$ intersects
  $\image(\pi_1(D,\nu)\to\pi_1(N,\nu))$ in a finitely generated
  subgroup, and
\item for each $i\in\{0,1\}$, and for each finitely generated subgroup
  $Z$ of $H_i$ which has the form $H_i\cap g^{-1}H_0g$ for some
  $g\in\pi_1(N,\nu)$, the based covering space of $(X_i,\nu)$
  corresponding to $Z$ admits a manifold compactification.
\EndConditions

Here conditions (ii) and (iii) hold trivially because $\pi_1(D)$ is
trivial. Condition (i) follows from the irreducibility of $N$. (A ball
bounded by a sphere in $\inter X_i$ must be contained in $X_i$ because
the frontier of $X_i$ is the disk $D$, and $\partial D\ne\emptyset$.)

To prove (iv), we consider any $i\in\{0,1\}$ and any subgroup of $H_i$
having the form $Z=H_i\cap g^{-1}H_0g$ where $g\in\pi_1(N,\nu)$.
Since $\pi_1(N,\nu)=H_0\star H_1$, we have either (a) $Z= \{1\}$ or
(b) $i=0$ and $g\in H_0$. If (a) holds then the based covering of
$(X_i,\nu)$ corresponding to $Z$ is equivalent to the universal cover
of $X_i$. But since $X_i$ is irreducible and has a non-empty boundary,
it is a Haken manifold. Hence by \cite[Theorem 8.1]{waldhausen}, the
universal cover of $X_i$ admits a manifold compactification. If (b)
holds then the covering corresponding to $Z$ is homeomorphic to $X_i$
and is therefore a manifold compactification of itself.
\EndProof

\Lemma\label{Simon consequence}
Suppose that $N$ is a compact, connected, orientable, irreducible PL
$3$-manifold and that $D$ is a separating, properly embedded disk in
$N$. Let $X$ denote the closure of one of the connected components of
$N-D$. Let $\nu\in D$ be a base point, and let $p:(\widetilde
N,\widetilde\nu)\to (N,\nu)$ denote the based covering space
corresponding to the subgroup $\image(\pi_1(X,\nu)\to\pi_1(N,\nu))$ of
$\pi_1(N,\nu)$. Let $\widetilde X$ denote the component of $p^{-1}(X)$
containing $\widetilde\nu$ (so that $p$ maps $\widetilde X$
homeomorphically onto $X$). Then every compact PL subset of
$\inter\widetilde N$ is PL ambient-isotopic to a subset of $\widetilde
X$.
\EndLemma

\Proof
Since $N$ is a compact, orientable, irreducible $3$-manifold
with non-empty boundary, it is a Haken manifold. Hence by
\cite[Theorem 8.1]{waldhausen}, the universal cover of $\inter N$ is
homeomorphic to ${\bf R}^3$. Thus $\inter\widetilde N$ is covered by
an irreducible manifold and is therefore irreducible.

According to Lemma \ref{Simon}, the manifold $\widetilde N$ admits a 
manifold compactification. Thus there is a homeomorphism $h$ of
$\widetilde N$ onto an open subset of a compact $3$-manifold $Q$ such
that $h(\inter\widetilde N)=\inter Q$. Since $\inter Q$ is
homeomorphic to the irreducible manifold $\inter\widetilde N$, the
compact manifold $Q$ is itself irreducible.
  
The definition of $\widetilde N$ implies that the inclusion map
$\iota:X\to N$ admits a based lift
$\widetilde\iota:(X,\nu)\to(\widetilde N,\widetilde\nu)$, that
$\widetilde\iota(X)=\widetilde X$, and that
$\widetilde\iota_\sharp:\pi_1(X,\nu)\to\pi_1(\widetilde
N,\widetilde\nu)$ is an isomorphism. Hence the inclusion $\widetilde
X\to \widetilde N$ induces an isomorphism of fundamental groups, and
if we set $X'=h(\widetilde X)$, the inclusion $X'\to Q$ induces an
isomorphism of fundamental groups.

On the other hand, since the frontier of $X$ in $N$ is $D$, the
frontier of {$X'$} in $Q$ is $D'= h(\widetilde\iota(D))$, a
properly embedded disk in the compact $3$-manifold $Q$. Set
$Y=\overline{Q-X'}$. Then in terms of a base point in $D'$ we have a
canonical identification of $\pi_1(Q)$ with
$\pi_1(X')\star\pi_1(Y)$. Since the inclusion $ X'\to Q$ induces an
isomorphism of fundamental groups, it follows that $\pi_1(Y)$ is
trivial. We also know that $Y$ is irreducible because its frontier in
the irreducible manifold $Q$ is a disk. Thus $Y$ is a compact, simply
connected, irreducible $3$-manifold with non-empty boundary, and is
therefore PL homeomorphic to a ball.

We have now exhibited $Q$ as the union of the compact $3$-dimensional
submanifold $X$ and the PL $3$-ball $Y$, and their intersection is the
disk $D$. It follows that any compact PL subset $W$ of $\inter Q$ is
PL isotopic to a subset of $\inter X$. Since $h$ maps
$\inter\widetilde N$ homeomorphically onto $\inter Q$, and maps
$\inter\widetilde X$ homeomorphically onto $\inter X'$, the conclusion
of the lemma follows.
\EndProof

\Lemma\label{handle removal lemma} Suppose that $K$ is a compact,
connected PL space such that $\pi_1(K)$ has rank $\ge2$ and is freely
indecomposable.  Suppose that $N$ is a compact, connected, orientable
PL $3$-manifold which is irreducible but boundary-reducible. Suppose
that $f:K\to\inter N$ is a $\pi_1$-injective PL map, and that the
inclusion homomorphism $\pi_1(f(K))\to\pi_1(N)$ is surjective. Then
$f$ is homotopic to a map $g$ such that $\ADS(g)<\ADS(f)$.  \EndLemma

\Proof
Since $N$ is boundary-reducible it contains an essential
properly embedded disk. If $N$ contains a non-separating essential
disk $D_0$, then there is a separating essential disk $D_1$ in $N-D_0$
such that the closure of the component of $N-D_1$ containing $D_0$ is
a solid torus $J$. In this case $\pi_1(\overline{N-J})$ is
non-trivial, since $\pi_1(N)$ has rank at least $2$; hence $D_1$ is an
essential disk as well. Thus in all cases, $N$ contains a separating
essential disk $D$. We may write $N=X_0\cup X_1$ for some connected
submanifolds $X_0$ and $X_1$ of $N$ with $X_0\cap X_1=D$. We choose a
base point in $\nu\in D$ and set
$A_i=\image(\pi_1(X_i,\nu)\to\pi_1(N,\nu))$ for $i=0,1$. Then
$\pi_1(N,\nu)=A_0\star A_1$.

If one of the $A_i$ were trivial, then one of the $X_i$ would be a
ball since $N$ is irreducible, and $D$ would not be an essential disk.
Hence the $A_i$ are non-trivial subgroups. It then follows from the
free product structure of $\pi_1(N,\nu)$ that the $A_i$ are of
infinite index in $\pi_1(N,\nu)$, and in particular that they are
proper subgroups.

Since the subgroup $H=f_\sharp(\pi_1(K))$, which is defined only up to
conjugacy in $\pi_1(N)$, has rank at least $2$ and is freely
indecomposable, it follows from the Kurosh subgroup theorem that $H$
is conjugate to a subgroup of one of the $A_i$. By symmetry we may
assume that $H$ is conjugate to a subgroup of $A_0$. Hence after
modifying $f$ by a homotopy we may assume that $f$ maps some base
point $\kappa$ of $K$ to $\nu$ and that
$f_\sharp(\pi_1(K,\kappa))\subset A_0$. Hence if $(\widetilde
N,\widetilde\nu)$ denotes the based covering space of $(N,\nu)$
corresponding to the subgroup $A_0$ of $\pi_1(N)$, then $f$ admits a
lift $\widetilde f:(K,\kappa)\to(\widetilde N,\widetilde\nu)$. Since
$A_0$ is a proper subgroup of $\pi_1(N,\nu)$, the covering space
$\widetilde N$ is non-trivial. Hence, according to Proposition
\ref{stallings}, we have $\ADS(\widetilde f)<\ADS(f)$.

Let $\widetilde X_0$ denote the component of $p^{-1}(X_0)$ containing
$\widetilde\nu$, so that $p$ maps $\widetilde X_0$ homeomorphically
onto $X_0$. According to Lemma \ref{Simon consequence}, the compact PL
subset $\widetilde f(K)$ of $\inter\widetilde N$ is PL
ambient-isotopic to a subset of $\widetilde X_0$.  In particular,
there is a PL homeomorphism $j$ of $\widetilde f(K)$ onto a subset $L$
of $\widetilde X_0$ such that $j$, regarded as a map of $\widetilde
f(K)$ into $\widetilde N$, is homotopic to the inclusion $\widetilde
f(K)\to\widetilde N$. It now follows that $p\circ j$ maps $\widetilde
f(K)$ homeomorphically onto the subset $p(L)$ of $X_0\subset N$. Hence
by Lemma \ref{homeomorphic images}, if we set $g=p\circ
j\circ\widetilde f:K\to N$, we have $\ADS(g)=\ADS(\widetilde
f)<\ADS(f)$. But since $j:\widetilde f(K)\to \widetilde N$ is
homotopic to the inclusion $\widetilde f(K)\to\widetilde N$, the map
$g:K\to N$ is homotopic to $f$.
\EndProof

\Proposition\label{DS and incompressible boundary} Suppose that $K$ is
a compact, connected PL space such that $\pi_1(K)$ has rank at least
$2$ and is freely indecomposable. Suppose that $f$ is a
$\pi_1$-injective.  PL map from $K$ to the interior of a compact,
connected, orientable, irreducible PL $3$-manifold $M$.  Then there
exist a map $g:K\to M$ homotopic to $f$ with $\ADS(g)\le\ADS(f)$, and
a compact, connected $3$-dimensional submanifold $N$ of $\inter M$
such that (i) $\inter N\supset g(K)$, (ii) the inclusion homomorphism
$\pi_1(g(K))\to\pi_1(N)$ is surjective, (iii) $\partial N$ is
incompressible in $M$, and (iv) $N$ is irreducible.  \EndProposition

\Proof
Among all maps from $K$ to $M$ that are homotopic to $f$, we
choose one, $g$, for which $\ADS(g)$ has the smallest possible value.
In particular we then have $\ADS(g)\le\ADS(f)$. Note also that
$f_\sharp:\pi_1(K)\to\pi_1(N)$ is injective.

Now let $N$ be a compact, connected $3$-submanifold of $M$ satisfying
conditions (i) and (ii) of the statement of the Proposition, and
choose $N$ so as to minimize the quantity $\kappa(\partial N)$
(see \ref{complexity def}) among all compact, connected
$3$-submanifolds satisfying (i) and (ii).   We shall
complete the proof by showing that $N$ satisfies (iii) and (iv).

We first show that (iv) holds, i.e. that $N$ is irreducible. If
$S\subset\inter N$ is a $2$-sphere, then $S$ bounds a ball $B\subset
M$. If we set $N'=N\cup B$, then the pair $N'$ satisfies (i) and
(ii). (It inherits property (ii) from $N$ because the inclusion
homomorphism $\pi_1(N)\to\pi_1(N')$ is surjective.) But if
$B\not\subset N$, it is clear from Definition \ref{complexity def}
that $\kappa(\partial N')<\kappa(\partial N)$, and the minimality of
$\kappa(\partial N)$ is contradicted. Hence we must have $B\subset N$,
and irreducibility is proved.

It remains to show that (iii) holds, i.e. that $\partial N$ is
incompressible. If this is false, then either $\partial N$ has a
sphere component, or there is a compressing disk $D$ for $\partial
N$. If $\partial N$ has a sphere component $S$, then the
irreducibility of $N$ implies that $N$ is a ball. But then the
injectivity of $g_\sharp:\pi_1(K)\to\pi_1(N)$ implies that $\pi_1(K)$
is trivial, a contradiction to the hypothesis that $\pi_1(K)$ has rank
at least $2$.
 
If there is a compressing disk $D$ for $\partial N$, then either
$D\cap N=\partial D$ or $D\subset N$. If $D\cap N=\partial D$, and if
we set $N'=N\cup Q$, where $Q$ is a regular neighborhood of $D$
relative to $\overline{M-N}$, then the { $3$-submanifold $N$ satisfies
conditions (i) and (ii).} (It inherits property (ii) because the
inclusion homomorphism $\pi_1(N)\to\pi_1(N')$ is again surjective.)
Now $\partial N'$ has the form $((\partial N)-A)\cup D_1\cup
D_2$, where $A\subset\partial N$ is a homotopically non-trivial
annulus, and $D_1$ and $D_2$ are disjoint disks in $M$ such that
$(D_1\cup D_2)\cap\partial N=\partial A$. It therefore follows from
Lemma \ref{pere tranquille} that $\kappa( \partial N)<\kappa(\partial
N')$. Again the minimality of $\kappa(\partial N)$ is contradicted.

Finally, if $D\subset N$, then $N$ is boundary-reducible. As we have
already shown that $N$ is irreducible, it follows from Lemma
\ref{handle removal lemma} that $f$ is homotopic in $N$ to a map $g'$
such that $\ADS(g')<\ADS(g)$. In particular, $g'$ is homotopic to $g$
in $M$; and since, according to \ref{co-restriction}, the absolute
degrees of singularity of $g$ and $g'$ do not depend on whether they
are regarded as maps into $N$ or into { $M$}, we now have a contradiction
to the minimality of $\ADS(g)$.
\EndProof

\section{Homology of covering spaces}
\label{homology section}

In this short section we shall apply and extend some  results from
\cite{shalenwagreich} concerning homology of covering spaces of
$3$-manifolds. In this section all homology groups are understood to
be defined with coefficients in $\Z_2$.

\Number\label{exact sequence}
If $N$ is a normal subgroup of a group $G$, we shall denote by $G\#N$
the subgroup of $G$ generated by all elements of the form
$gag^{-1}a^{-1}b^2$ with $g \in G$ and $a,b \in N$. (This is a special
case of the notation used in \cite{Stallingshomology} and \cite
{shalenwagreich}. Here we are taking the prime $p$, which was
arbitrary in \cite{Stallingshomology} and \cite {shalenwagreich}, to
be $2$.
\EndNumber

\Number\label{gee sub em} As in Section 1 of \cite{shalenwagreich},
for any group $\Gamma$, we define subgroups $\Gamma_{d}$ of
$\Gamma$ recursively for ${d}\ge0$, by setting
$\Gamma_0=\Gamma$ and
$\Gamma_{{d}+1}=\Gamma\#\Gamma_{d}$. We regard
$\Gamma_{d}/\Gamma_{{d}+1}$ as a $\Z_2$-vector space.
\EndNumber

\Lemma \label{2r-3}
Let $M$ be a closed $3$-manifold and set $r=\rk M$.  Suppose that
$\widetilde M$ is a regular cover of $M$ whose group of deck
transformations is isomorphic to $(\Z_2)^m$ for some integer $m\ge0$.
Then
$$\rk(\widetilde M)\ge mr-\frac{m(m+1)}2.$$
\EndLemma

\Proof
We set $\Gamma=\pi_1(M)$ and define $\Gamma_{d}$ for each ${d}\ge0$ as in
\ref{exact sequence}. We have $\rk \Gamma/\Gamma_1=\rk M=r$. It then
follows from \cite[Lemma 1.3]{shalenwagreich} that
$\rk(\Gamma_1/\Gamma_2) \ge r(r-1)/2$.

Let $N$ denote the normal subgroup of $\Gamma$ corresponding to the
regular covering space $\widetilde M$. Since
$\Gamma/N\cong(\Z_2)^m$, we may write $N = E\Gamma_1$ for some
$(r-m)$-generator subgroup $E$ of $\Gamma$.  It now follows from
\cite[Lemma 1.4]{shalenwagreich} that
\begin{eqnarray*}
\rk\widetilde M&=&\rk H_1(E\Gamma_1)\cr
&\ge&\rk(\Gamma_1/\Gamma_2)-\frac{(r-m)(r-m-1)}{2}\cr
&\ge &
\frac{r(r-1)}{2}-\frac{(r-m)(r-m-1)}2 = mr-\frac{m(m+1)}2.
\end{eqnarray*}
\EndProof

The case $m=2$ of Lemma \ref{2r-3} will be applied in the proof of
Lemma \ref{2r-4}.

\section{An application of a result of Gabai's}
\label{new section}

This section contains the applications of Gabai's results that were
mentioned in the introduction.  The main result of the section
is Proposition \ref{newprop}.

\Lemma\label{dark green}
Let $X$ be a PL space, let $K$ be a closed, connected, orientable
surface of genus $g>0$, and let $f:K\to X$ be a PL map. Suppose that
the homomorphism
$f_*:H_2(K;\Z_2)\to H_2(X;\Z_2)$ is trivial. Then the image of
$f_*:H_1(K;\Z_2)\to H_1(X;\Z_2)$ has dimension at most
$g$.
\EndLemma

\Proof
Since $f_*:H_2(K;\Z_2)\to H_2(X;\Z_2)$ is trivial, it follows that the dual
homomorphism $f^*:H^2(X;\Z_2)\to H^2(K;\Z_2)$ is also trivial. Hence
for any $\alpha,\beta\in H^1(X)$ we have
$$f^*(\alpha)\cup f^*(\beta)=f^*(\alpha\cup\beta)=0.$$
Thus if we set $V=H^1(K;\Z_2)$ and let $L\subset V$ denote the image of
of $f^*:H^1(X;\Z_2)\to H^1(K;\Z_2)$, we have $L\cup L=0$, i.e.
$$L\subset L^\perp=\{v\in V:v\cup L=0\}.$$ 
Hence if $d$ denotes the dimension of $L$, we have 
$$d\le\dim L^\perp.$$
But by Poincar\'e duality, the cup product pairing on $V$ is
non-singular, and so
$$\dim L^\perp=\dim V-\dim L=2g-d.$$
Hence $d\le g$. 
As the linear map $f_*:H_1(K;\Z_2)\to H_1(X;\Z_2)$ is dual to 
$f^*:H^1(X;\Z_2)\to H^1(K;\Z_2)$, its rank is the same as that of
$f^*$, namely $d$. The conclusion follows.
\EndProof

\Notation
If $F$ is a closed, orientable $2$-manifold, we shall denote by $\tg(F)$ the
total genus of $F$, i.e. the sum of the genera of its
components.
\EndNotation

\Lemma\label{polka dots}
For any compact, connected, orientable $3$-manifold  $N$, we have
$$\tg(\partial  N)\le\rk  N.$$
\EndLemma

\Proof
In the exact sequence
$$H_2( N,\partial N;\Z_2)\longrightarrow H_1(\partial N;\Z_2)
\longrightarrow H_1( N;\Z_2) ,$$
Poincar\'e-Lefschetz duality implies that the vector spaces
$H_2(N,\partial N;\Z_2)$ and $H_1( N;\Z_2)$ are of the same dimension,
$\rk N$.  Hence we have
$$2~\tg(\partial  N)=\rk\partial  N\le2\rk  N$$
and the conclusion follows.
\EndProof

\Lemma\label{falafel} If $N$ is a compact, connected, orientable
$3$-manifold $N$ such that $\partial N$ has at most one connected
component, then $H_2(N;\Z)$ is torsion-free.  \EndLemma

\Proof
In the exact sequence
$$H_2(\partial N;\Z)\longrightarrow H_2( N;\Z)\longrightarrow H_2(
N,\partial N;\Z)$$
the inclusion map $H_2(\partial
N;\Z)\longrightarrow H_2( N;\Z)$ is trivial since $\partial N$ has at
most one connected component. Hence the map $H_2( N;\Z)\longrightarrow
H_2( N,\partial N;\Z)$ is injective, so that $H_2( N;\Z)$ is
isomorphic to a subgroup of $ H_2( N,\partial N;\Z)$. But by
Poincar\'e-Lefschetz duality, $ H_2( N,\partial N;\Z)$ is isomorphic
to $H^1(N,\Z)$ and is therefore torsion-free. The conclusion follows.
\EndProof

\Proposition\label{newprop}
Suppose that $N$ is a compact (possibly closed) orientable
$3$-manifold which is irreducible and boundary-irreducible.  Suppose
that $K$ is a closed, connected, orientable surface of genus $g\ge2$,
and that $\phi:K\to N$ is a $\pi_1$-injective PL map.  Then either
\Conclusions
\item $N$ contains a connected (non-empty) closed incompressible
  surface of genus at most $g$, or
\item the $\Z_2$-vector subspace $\phi_*(H_1(K;\Z_2))$ of
  $H_1(N;\Z_2)$ has dimension at most $g$.
\EndConclusions
Furthermore, if $\phi_*:H_1(K;\Z_2)\to H_1(N;\Z_2)$ is surjective and
$\partial N\ne\emptyset$, then (i) holds.  \EndProposition

\Proof
We begin with the observation that $N$ is non-simply connected in view
of the existence of the map $\phi$. Since $N$ is also irreducible, it
follows that no component of $\partial N$ is a sphere.
On the other hand, since $N$ is boundary-irreducible, every component
of $\partial N$ is $\pi_1$-injective in $N$. Thus every component
of $\partial N$ is parallel to an incompressible surface in $N$. 

To prove the first assertion of the proposition we distinguish three
cases, which are not mutually exclusive but cover all possibilities.
\Cases 
\Case{{\bf Case A.\ }} The homomorphism $\phi_*:H_2(K;\Z)\to
H_2(N;\Z)$ is trivial.
\Case{{\bf Case B.\ }} The surface $\partial N$
has at least two components.
\Case{{\bf Case C.\ }} The surface
$\partial N$ has at most one component and
$\phi_*:H_2(K;\Z)\to H_2(N;\Z)$ is a non-trivial homomorphism.
\EndCases

To prove the assertion in Case A, we first consider the commutative diagram
$$\xymatrix{
H_2(K;\Z)    \ar[r]\ar[d]  &     H_2(N;\Z) \ar[d] \\
H_2(K;\Z_2)  \ar[r]        &     H_2(N;\Z_2) \\
}$$
in which the vertical maps are natural homomorphisms and the
horizontal maps are induced by $\phi$. The left-hand vertical arrow is
surjective because the surface $K$ is orientable.  Since the top
horizontal map is trivial, it follows that the bottom horizontal map
is trivial. Hence Lemma \ref{dark green} asserts that the image of
$\phi_*:H_1(K;\Z_2)\to H_1(N;\Z_2)$ has dimension at most $g$. Thus
alternative (2) of the conclusion holds in Case A.

In Case B, using Lemma \ref{polka dots} and the surjectivity of
$\phi_*:H_1(K;\Z_2)\to H_1(N;\Z_2)$, we find that
$$\tg(\partial N)\le \rk N\le\rk K=2g.$$
Since $\partial N$ has at
least two components in this case, some component $F$ of $\partial N$
must have genus at most $g$.  By the observation at the beginning of the
proof, $F$ is parallel to an incompressible surface in $N$. Thus
alternative (1) of the conclusion holds in Case B.

To prove the assertion in Case C, we begin by considering the commutative diagram
$$\xymatrix{
H_2(K;\Z)  \ar[r]\ar[d]   &   H_2(N;\Z) \ar[d] \\
H_2(K;\R)  \ar[r]         &   H_2(N;\R) \\
}$$
in which the vertical maps are natural homomorphisms and the
horizontal maps are induced by $\phi$. Since $\partial M$
has at most one component, Lemma \ref{falafel} asserts that $H_2(N;\Z)$ is
torsion-free. Hence the right-hand vertical arrow in the diagram is
injective.  Since the top horizontal map is non-trivial, it follows
that the bottom horizontal map is non-trivial. In other words, if
$[K]$ denotes the fundamental class in $H_2(K;\R)$ then the class
$\alpha=f_*([K])\in H_2(N;\R)$ is non-zero.

We shall now apply a result from \cite{gabai}.  For any $2$-manifold
$\mathcal F$ we shall denote by $\chiminus({\mathcal F})$ the quantity 
$$\sum_F\max(\chibar(F),0),$$
where $F$ ranges over the components of
$\mathcal F$.  As in \cite{gabai}, given a class $z$ in $H_2(M;\R)$,
we denote by $x_s(z)$ and $x(z)$ respectively the ``norm based on
singular surfaces'' and the Thurston norm of $z$.  Since $\alpha$ is
by definition realized by a map of the surface $K$ into $N$, and since
$\chiminus(K)=2g-2$, we have $x_s(\alpha)\le 2g-2$.  But it follows
from \cite[Corollary 6.18]{gabai} that $x(\alpha) = x_s(\alpha)$.
Hence $x(\alpha)\le 2g-2$.  By definition this means that if
${\mathcal F}$ is a closed orientable embedded surface in $\inter N$
such that the fundamental class $[{\mathcal F}]\in H_1({\mathcal
  F};\R)$ is mapped to $\alpha$ under inclusion, and if $\mathcal F$
is chosen among all such surfaces so as to minimize
$\chiminus({\mathcal F})$, then $\chiminus({\mathcal F})\le 2g-2$.
Since $\alpha\ne0$ we have ${\mathcal F}\ne\emptyset$.

Since $N$ is  irreducible, any sphere component of
$\mathcal F$ must be homologically trivial in $N$. We may assume
that every torus component of $F$ is compressible, as otherwise
alternative (1) of the conclusion holds. Under this assumption, if $T$ is a torus
component of $\mathcal F$, compressing
$T$ yields a sphere which must be homologically trivial; hence $T$ is
itself homologically trivial. Thus after discarding homologically
trivial components of $\mathcal F$ whose Euler characteristics are
$\ge0$, we may suppose that no component of $\mathcal F$ is a sphere or
torus. The minimality of $\chiminus({\mathcal F})$ now implies that
$\mathcal F$ is incompressible.

 Let $F$ be any
component of $\mathcal F$. Then $F$ is an incompressible closed surface in
$N$, and we have
$$\chiminus(F)\le\chiminus({\mathcal F})\le2g-2.$$
Hence $F$ has genus at most $g$, and alternative (1) holds. This
completes the proof of the first assertion of the proposition.

To prove the second assertion, suppose that
 $\phi_*:H_1(K;\Z_2)\to
H_1(N;\Z_2)$ is  surjective, that $\partial N\ne\emptyset$, and that
alternative (2) holds. Then $\rk N\le
g$, and it follows from Lemma \ref{polka dots} that
$\tg(\partial N)\le g$. In particular, any component $F$ of the
non-empty $2$-manifold $\partial N$ has genus at most $g$. 
By the observation at the beginning of the
proof, $F$ is parallel to an incompressible surface in $N$. Thus
alternative (1) of the conclusion holds.
\EndProof

\section{Towers}
\label{tower section}

In this section we prove a result, Proposition \ref{there's an
  extension}, which summarizes the tower construction described in the
introduction. Our main topological result, Theorem \ref{top 11}, will
then be proved by combining Proposition \ref{there's an extension}
with results from the earlier sections.  We begin by introducing some
machinery that will be needed for the statement and proof of
Proposition \ref{there's an extension}.

\Definition\label{tower def}
Suppose that $n$ is a non-negative integer.  We define a {\it
height-$n$ tower of $3$-manifolds} to be a $(3n+2)$-tuple
$${\mathcal T}=(M_0,N_0,p_1,M_1,N_1,p_2,\ldots,p_n,M_n,N_n),$$ where
$M_0,\ldots,M_n$ are compact, connected, orientable PL $3$-manifolds,
$N_j$ is a compact, connected $3$-dimensional PL submanifold of $M_j$
for $j=0,\ldots,n$, and $p_j:M_j\to N_{j-1}$ is a PL covering map for
$j=1,\ldots,n$.  We shall refer to $M_0$ as the {\it base} of the
tower $\mathcal T$ and to $N_n$ as its {\it top}.  We define the {\it
tower map associated to $\mathcal T$} to be the map
$$h=\iota_0\circ p_1\circ\iota_1\circ p_2\circ\cdots\circ
p_n\circ\iota_n:N_n\to M_0,$$
where $\iota_j:N_j\to M_j$ denotes the inclusion map  for
$j=0,\ldots,n$.
\EndDefinition

\Number\label{tower remark}
Consider any tower of $3$-manifolds
$${\mathcal T}=(M_0,N_0,p_1,M_1,N_1,p_2,\ldots,p_n,M_n,N_n).$$
Note that for any given $j$ with $0\le
j<n$, the manifold $N_j$ is closed if and only if its finite-sheeted
covering space $M_{j+1}$ is closed. Note also that if, for a given $j$
with $0\le j\le n$, the submanifold $N_j$ of the (connected) manifold
$M_j$ is closed, then we must have $N_j=M_j$, so that in particular
$M_j$ is closed.  

It follows that if $M_j$ is closed for a given $j$ with $0\le j\le n$,
then $M_{i}$ is also closed for every $i$ with $0\le i\le j$. Thus
either all the $M_j$ { have non-empty {boundaries}}, or there is an
index $j_0$ with $0\le j_0\le n$ such that $M_j$ is closed when $0\le
j\le j_0$ and {$M_j$ has non-empty boundary} when $j_0<j\le n$.
Furthermore, in the latter case, for each $j< j_0$ we have $N_j=M_j$.
\EndNumber

\Number\label{further tower remark}
In particular, if in a tower of $3$-manifolds
$${\mathcal T}=(M_0,N_0,p_1,M_1,N_1,p_2,\ldots,p_n,M_n,N_n)$$ the
manifold $M_j$ is closed for a given $j\le n$, then for every $i$ with
$0\le i<j$ the composition
$$p_{j-1}\circ\cdots\circ p_{i}:M_j\to M_i$$ is a well-defined
covering map, whose degree is the product of the degrees of
$p_i,\ldots,p_{j-1}$.
\EndNumber

\Definition\label{good tower def}
A tower of $3$-manifolds
$${\mathcal T}=(M_0,N_0,p_1,M_1,N_1,p_2,\ldots,p_n,M_n,N_n)$$
 will be termed {\it good} if it has the following 
properties:
\Properties
\item $M_j$ and $N_j$ are irreducible for $j=0,\ldots,n$;
\item $\partial N_j$ is a (possibly empty) incompressible surface in
  $M_j$ for $j=0,\ldots,n$;
\item the covering map $p_j:M_j\to N_{j-1}$ has degree $2$ for
  $j=1,\ldots,n$; and
\item for each $j$ with $2\le j\le n$ such that $M_j$ is closed, the
  four-fold covering map (see \ref{further tower remark})
  $$p_j\circ p_{j-1}:M_j\to M_{j-2}$$
  is regular and has covering group isomorphic to $\Z_2\times\Z_2$.
\EndProperties
\EndDefinition

\Lemma\label{2r-4}
Suppose that
$${\mathcal T}=(M_0,N_0,p_1,M_1,N_1,p_2,\ldots,p_n,M_n,N_n)$$
is a good tower of $3$-manifolds and that $j_0$ is an index with
$0 \le j_0\le n$ such that $M_{j_0}$ is closed.
Set $r=\rk M_0$ and assume that $r\ge3$.
For any index  $j$ with $0\le j\le j_0$, we have 
$$\rk M_j\ge2^{j/2}(r-3)+3$$
if $j$ is even, and
$$\rk M_j\ge2^{(j-1)/2}(r-3)+2$$
if $j$ is odd. 

In particular, we have $\rk M_j\ge r-1$ for each $j$ with $0\le j\le
n$ such that $M_j$ is closed, and we have $\rk M_j\ge 2r-4$ for each
$j$ with $2\le j\le n$ such that $M_j$ is closed.
\EndLemma

\Proof
According to \ref{tower remark}, $M_j$ is closed for every index $j$
with $0\le j \le j_0$.  We shall first show that for every even $j$
with $0\le j\le j_0$ we have $\rk M_j\ge2^{j/2}(r-3)+3$. For $j=0$
this is trivial since $r=\rk M_0$. Now, arguing inductively, suppose
that $j$ is even, that $0< j\le n$, and that $\rk
M_{j-2}\ge2^{(j-2)/2}(r-3)+3$.  Since the definition of a good tower
implies that $M_j$ is a regular $(\Z_2\times\Z_2)$-cover of $M_{j-2}$,
we apply Lemma \ref{2r-3} with $m=2$ to deduce that
$$\rk M_j\ge2(\rk M_{j-2})-3\ge
2(2^{(j-2)/2}(r-3)+3)-3=2^{j/2}(r-3)+3.$$ This completes the induction
and shows that $\rk M_j\ge2^{j/2}(r-3)+3$ for every even index $j$
with $2\le j\le j_0$.  Finally, if $j$ is an odd index with $2<j\le
j_0$, then since $j-1$ is even we have $\rk
M_{j-1}\ge2^{(j-1)/2}(r-3)+3$; and since $M_j$ is a $2$-sheeted cover
of $M_{j-1}$, it is clear that $\rk M_j\ge\rk
M_{j-1}-1\ge2^{(j-1)/2}(r-3)+2$.
\EndProof

\Definition\label{truncation def}
If 
$${\mathcal T}=(M_0,N_0,p_1,M_1,N_1,p_2,\ldots,p_n,M_n,N_n)$$ is a
height-$n$ tower of $3$-manifolds, then for any $m$ with $0\le
m\le n$, the $(3m+2)$-tuple
$${\mathcal T}^-=(M_0,N_0,p_1,M_1,N_1,p_2,\ldots,p_m,M_m,N_m)$$
is a
height-${n}$ tower. We shall refer to the tower ${\mathcal T}^-$ as the
height-$m$ {\it truncation} of $\mathcal T$. We shall say that a tower
${\mathcal T}^+$ is an {\it extension} of a tower $\mathcal T$, or
that 
${\mathcal T}^+$ {\it extends} $\mathcal T$, if ${\mathcal T}$ is a
truncation of ${\mathcal T}^+$.

In particular, any tower may be regarded as an extension of
  itself. This will be called the {\it degenerate} extension.
\EndDefinition

\Definition\label{good lift def}
Let $\mathcal T$ be a tower of $3$-manifolds with base $M$ and
top $N$, and let $h:N\to M$ denote the associated tower map. Let $\phi$
be a PL map from a compact PL space $K$ to $M$. By a {\it homotopy-lift}
of $\phi$ through the tower $\mathcal T$ we mean a PL map $\widetilde \phi:K\to
N$ such that $h\circ\widetilde \phi$ is homotopic to $\phi$. A homotopy-lift
$\widetilde \phi$ of $\phi$ will be termed {\it good} if the inclusion
homomorphism $\pi_1({\widetilde\phi(K)})\to\pi_1(N)$ is surjective.
\EndDefinition

\Lemma\label{first hairy lemma}
Suppose that $K$ is a compact PL space with freely indecomposable
fundamental group of rank $k\ge2$.  Suppose that $\mathcal
T=(M_0,N_0,p_1,\ldots,N_n)$ is a good tower of $3$-manifolds of height
$n$.  Suppose that $\phi:K\to M_0$ is a $\pi_1$-injective PL map, and that
$\widetilde \phi:K\to N_n$ is a good homotopy-lift of $\phi$ through
the tower $\mathcal T$. Suppose that $p_{n+1}:M_{n+1}\to N_n$ is a
two-sheeted covering space of $N_n$, and that the map $\widetilde
\phi:K\to N_n$ admits a lift to the covering space $M_{n+1}$. Suppose
that {\it either}
\begin{enumerate}
\item[($\alpha$)]$n\ge1$, the manifold $N_{n}$ is closed (so that
$M_{n+1}$ is closed, cf.  \ref{tower remark}), and the covering map
$$p_n\circ p_{n+1}:M_{n+1}\to M_{n-1}$$ is regular and has covering
group isomorphic to $\Z_2\times\Z_2$; or
\item[($\beta$)]$\partial M_{n+1}\ne\emptyset$,  or
\item[($\gamma$)]$n=0$.
\end{enumerate}
Then there exists a compact submanifold $N_{n+1}$ of $M_{n+1}$ with
the following properties: 
\Conclusions
\item ${\mathcal
    T}^+=(M_0,N_0,p_1,\ldots,N_n,p_{n+1},M_{n+1},N_{n+1})$ is a good
  height-$(n+1)$ tower extending ${\mathcal T}$, and
\item there is a a good homotopy-lift $\widetilde \phi^+$ of $\phi$
  through the tower ${\mathcal T}^+$ such that
$$\ADS(\widetilde \phi^+)<\ADS(\widetilde \phi).$$
\EndConclusions
\EndLemma

\Proof
Let $h:N_n\to M_0$ be the tower map associated to $\calt$. We fix a
lift $f:K\to M_{n+1}$ of the map $\widetilde \phi:K\to N_n$ to the
covering space $M_{n+1}$. Since $\widetilde\phi$ is a homotopy lift of
$\phi$, the map $h\circ p_{n+1}\circ f:K\to M_0$ is homotopic to
$\phi$.  Since $\phi_\sharp:\pi_1(K)\to\pi_1( M_0)$ is injective, it
now follows that $f_\sharp:\pi_1(K)\to\pi_1(M_{n+1})$ is also
injective.  We may therefore apply Proposition \ref{DS and
incompressible boundary} to this map $f$, taking $M=M_{n+1}$ and
$N=N_{n+1}$. We choose a map $g:K\to M_{n+1}$ homotopic to $f$, with
$\ADS(g)\le\ADS(f)$, and a compact $3$-dimensional submanifold
$N=N_{n+1}$ of $\inter M_{n+1}$, such that conditions (i)---(iv) of
\ref{DS and incompressible boundary} hold with $M=M_{n+1}$.

It is clear from the definition that ${\mathcal
T}^+=(M_0,N_0,p_1,\ldots,N_n,p_{n+1},M_{n+1},N_{n+1})$ is a tower
extending $\calt$. To show that the tower $\calt^+$ is good, we first
observe that conditions (1)---(4) of Definition \ref{good tower def}
hold whenever $j\le n$ because $\calt$ is a good tower. For $j=n+1$,
Conditions (1) and (2) of Definition \ref{good tower def} follow from
conditions (iv) and (iii) of \ref{DS and incompressible boundary},
while condition (3) of Definition \ref{good tower def} follows from
the hypothesis that $p_{n+1}:M_{n+1}\to N_n$ is a two-sheeted
covering. The case $j=n+1$ of Condition (4) of Definition \ref{good
tower def} is clear if alternative ($\alpha$) of the hypothesis holds,
and is vacuously true if alternative ($\beta$) or ($\gamma$)
holds. Hence $\calt^+$ is a good tower.

Since by condition (i) of Proposition \ref{DS and incompressible
boundary} we have $\inter N_{n+1}\supset g(K)$, we may regard
$g:K\to M_{n+1}$ as a composition $\iota_{n+1}\circ\widetilde \phi^+$,
where $\iota_{n+1}:N_{n+1}\to M_{n+1}$ is the inclusion map and
$\widetilde \phi^+$ is a PL map from $K$ to $N_{n+1}$.  Since $g$ is
homotopic to $f$, the map {$h\circ
p_{n+1}\circ\iota_{n+1}\circ\widetilde \phi^+= h\circ p_{n+1}\circ
g:K\to M_0$} is homotopic to $\phi$.  It follows that
$\widetilde\phi^+$ is a homotopy-lift of $\phi$ through the tower
$\calt^+$. Condition (ii) of \ref{DS and incompressible boundary}
asserts that the inclusion homomorphism
$\pi_1(\widetilde\phi^+(K))\to\pi_1(N_{n+1})$ is surjective, which
according to Definition \ref{good lift def} means that the
homotopy-lift $\widetilde\phi^+$ of $\phi$ is good.

Finally, since the homotopy-lift $\widetilde\phi$ of $\phi$ is good by
hypothesis, the inclusion homomorphism
$\pi_1(\widetilde\phi^+(K))\to\pi_1(N_n)$ is surjective. As $f$ is a
lift of $\widetilde\phi$ to the non-trivial covering space $M_{n+1}$
of $N_n$, it follows from Proposition \ref{stallings} that
$\ADS(f)<\ADS(\widetilde\phi)$. But we chose $g$ in such a way that
$\ADS(g)\le\ADS(f)$, and according to \ref{co-restriction} we have
$\ADS(\widetilde\phi^+)=\ADS(g)$. Hence
$\ADS(\widetilde\phi^+)<\ADS(\widetilde\phi)$.
\EndProof

\Lemma\label{second hairy lemma}
Suppose that $K$ is a closed orientable surface of genus
$g\ge2$.  Suppose that $\mathcal T$ is a good tower of $3$-manifolds
of height $n$.  Let $M$ denote the base of $\mathcal T$,
and assume that $\rk M\ge g+3$. Suppose
that $\phi:K\to M$ is a $\pi_1$-injective PL map, and that $\widetilde
\phi$ is a good homotopy-lift of $\phi$ through the tower $\mathcal
T$. Then at least one of the following alternatives holds:
\Conclusions
\item $N_n$ contains a connected (non-empty) closed incompressible
  surface of genus at most $g$;
\item $n\ge1$ and $N_{n-1}$ contains a connected (non-empty) closed
  incompressible surface of genus at most $g$; or
\item there exist a height-$(n+1)$ extension ${\mathcal T}^+$ of
  $\mathcal T$ which is a good tower, and a good homotopy-lift
  $\widetilde \phi^+$ of $\phi$ through the tower ${\mathcal T}^+$,
  such that
$$\ADS(\widetilde \phi^+)<\ADS(\widetilde \phi).$$
\EndConclusions
\EndLemma

\Proof
We write
$${\mathcal T}=(M_0,N_0,p_1,M_1,N_1,p_2,\ldots,p_n,M_n,N_n),$$ 
so that $M=M_0$. We distinguish several cases.

\Cases 
\Case{{\bf Case A:\ }}$\partial N_n\ne\emptyset$ and the
homomorphism $\widetilde \phi_*:H_1(K;\Z_2)\to H_1(N_n;\Z_2)$ is
surjective;
\Case{{\bf Case B:\ }} $\partial N_n\ne\emptyset$ and
$\widetilde \phi_*:H_1(K;\Z_2)\to H_1(N_n;\Z_2)$ is not surjective;
\Case{{\bf Case C:\ }}$n=0$;
\Case{{\bf Case D:\ }}$n\ge1$ and $N_n$ is closed.
\EndCases

In Case A, all the hypotheses of the final assertion of
Proposition \ref{newprop} hold with $\widetilde\phi$ in place of
$\phi$. It therefore follows from the final assertion of
Proposition \ref{newprop} that conclusion (1) of the present lemma holds.

In Case B, the map $\widetilde\phi:K\to N_n$ admits a lift to some
two-sheeted covering space $p_{n+1}:M_{n+1}\to N_n$ of $N_n$. Since
$\partial N_n\ne\emptyset$, we have $\partial M_{n+1}\ne\emptyset.$
This is alternative ($\beta$) of the hypothesis of Lemma \ref{first
hairy lemma}. It therefore follows from \ref{first hairy lemma} that
conclusion (3) of the present lemma holds.

In Case C the argument is identical to the one used in Case B, except
that we have alternative $(\gamma)$ of Lemma \ref{first hairy lemma}
in place of alternative ($\beta$).

We now turn to Case D. In this case, as was observed in \ref{tower
  remark}, we have $N_n=M_n$ and $N_{n-1}=M_{n-1}$, and $p_n$ is a
two-sheeted covering map from $M_n$ to $M_{n-1}$.

Let us set $r=\rk M\ge g+3$. According to Lemma \ref{2r-4}, for any
index $j$ such that $1\le j\le n$ and such that $M_j$ is closed, we
have $\rk M_j\ge r-1$. In particular, if we set $d=\rk M_{n-1}$, we
have $d\ge r-1\ge g+2$.

Now set $\widetilde\phi^-=p_n\circ\widetilde\phi :K\to M_{n-1}$.  Then
$X=\widetilde\phi^-_*(H_1(K;\Z_2))$ is a subspace of the
$d$-dimensional $\Z_2$-vector space $V=H_1(M_{n-1};\Z_2)$.

The hypotheses of Proposition \ref{newprop} hold with $N$ and
$\phi$ replaced by $M_{n-1}$ and $\widetilde\phi^-$. Hence either
$M_{n-1}$ contains a connected (non-empty) closed incompressible
surface of genus at most $g$, or $X$ has dimension
at most $g$. The first alternative gives conclusion (2) of the present
lemma.

There remains the subcase in which $X$ has dimension at most $g$.
Since $d\ge r-1\ge g+2$, the dimension of $X$ is then at most $g\le
d-2$.

In this subcase we shall show that $\widetilde\phi:K\to
M_n$ admits a lift to some two-sheeted covering space
$p_{n+1}:M_{n+1}\to M_n$ of $M_n=N_n$ for which alternative ($\alpha$)
of the hypothesis of Lemma \ref{first hairy lemma} holds. It will
 will then  follow from \ref{first hairy lemma} that
conclusion (3) of the present lemma holds.

Let $q$ denote the natural homomorphism from $\pi_1(M_{n-1})$ to
$H_1(M_{n-1};\Z_2)$. The two-sheeted cover $M_n$ of $M_{n-1}$
corresponds to a normal subgroup of $\pi_1(M_{n-1})$ having the form
$q^{-1}(Z)$, where $Z$ is some $(d-1)$-dimensional vector subspace of
$V$. Since $\widetilde\phi^-$ admits the lift $\widetilde\phi$ to
$M_n$, we have $X\subset Z\subset V$. Since in addition we have $\rk
X\le d-2<d-1=\rank Z$, there exists a $(d-2)$-dimensional vector
subspace $Y$ of $V$ with $X\subset Y\subset Z$. The subgroup
$q^{-1}(Y)$ determines a regular covering space $ M_{n+1}$ of $M_{n-1}$
with covering group $\Z_2\times\Z_2$. Since $q^{-1}(Y)\subset
q^{-1}(Z)$, the degree-four covering map $M_{n+1}\to M_{n-1}$ factors
as the composition of a degree-two covering map $p_{n+1}:M_{n+1}\to
M_n$ with $p_n:M_n\to M_{n+1}$. Thus the covering space
$p_{n+1}:M_{n+1}\to M_n$ satisfies alternative ($\alpha$) of
\ref{first hairy lemma}. It remains to show that $\widetilde\phi$
admits a lift to $M_{n+1}$.

Since $\widetilde\phi^-_\sharp(\pi_1(K))\subset q^{-1}(X)\subset
q^{-1}(Y)$, the map $\widetilde\phi^-$ admits a lift to the four-fold
cover $M_{n+1}$ of $M_{n-1}$. Since $M_{n+1}$ is a regular covering
space of $M_{n-1}$, there exist four different lifts of $\phi^-$ to
$M_{n+1}$. But $\widetilde\phi^-$ can have at most two lifts to $M_n$,
and each of these can have at most two lifts to $M_{n+1}$. Hence each
lift of $\widetilde\phi^-$ to $M_n$ admits a lift to $M_{n+1}$. In
particular, $\widetilde\phi$ admits a lift to $M_{n+1}$.
\EndProof

\Lemma\label{there's an extension} Suppose that $K$ is a closed,
orientable surface of genus $g\ge2$.  Suppose that $M$ is a closed,
orientable $3$-manifold such that $\rk M\ge g+3$, and that $\phi:K\to
M$ is a $\pi_1$-injective PL map. Suppose that
$$\calt_0=(M_0,N_0,p_1,\ldots,N_{n_0})$$
is a good tower with base $M$
such that $\phi$ admits a good homotopy-lift through $\calt$.  Then
either 
\Conclusions 
\item ${n_0}\ge1$, and $N_{{n_0}-1}$ contains a connected (non-empty)
  closed incompressible surface of genus at most $g$, or
\item there exists a good tower $\calt_1$ which is a (possibly
  degenerate) extension of $\calt_0$, such that the top $N$ of
  $\calt_1$ contains a connected (non-empty) closed incompressible
  surface of genus at most $g$, and $\phi$ admits a good homotopy-lift
  $\widetilde\phi_1$ through the tower $\calt_1$.
\EndConclusions
\EndLemma

\Proof 
Let us fix a good homotopy-lift $\widetilde\phi_0$ of $\phi$ through
$\calt_0$.  Let $\cals$ denote the set of all ordered pairs
$(\calt,\widetilde\phi)$ such that $\calt$ is a good tower which is an
extension of $\calt_0$ and $\widetilde\phi$ is a good homotopy-lift of
$\phi$ through $\calt$. Then we have $(\calt_0,\phi_0)\in\cals$, and
so $\cals\ne\emptyset$. Hence there is an element
$(\calt_1,\widetilde\phi_1)$ of $\cals$ such that
$\ADS(\widetilde\phi_1)\le\ADS(\widetilde\phi)$ for every element
$(\calt,\widetilde\phi)$ of $\cals$. Let us write
$$\calt_1=(M_0,N_0,p_1,\ldots,N_{n_1}).$$

The hypotheses of Lemma \ref{second hairy lemma} now hold with
  $\calt_1$ and $\widetilde\phi_1$ in place of $\calt$ and
  $\widetilde\phi$. Hence one of the following alternatives must hold:
\begin{enumerate}
\item[\ref{second hairy lemma}(1)] $N_{n_1}$ contains a connected
  (non-empty) closed incompressible surface of genus at most $g$;
\item[\ref{second hairy lemma}(2)]$n_1\ge1$ and $N_{n_1-1}$ contains a
  connected (non-empty) closed incompressible surface of genus at most
  $g$; or
\item[\ref{second hairy lemma}(3)] there exist a height-$(n_1+1)$
  extension ${\mathcal T}^+$ of $\calt_1$ which is a good tower, and a
  good homotopy-lift $\widetilde \phi^+$ of $\phi$ through the tower
  ${\mathcal T}^+$, such that
$$\ADS(\widetilde \phi^+)<\ADS(\widetilde \phi_1).$$
\end{enumerate}

If \ref{second hairy lemma}(1) holds, then the tower $\calt_1$ has the
property asserted in conclusion (1) of the present lemma. If
\ref{second hairy lemma}(2) holds, and if $n_1>{n_0}$ (i.e. $\calt_1$
is a non-degenerate extension of $\calt_0$), then the height-$(n_1-1)$
truncation $\calt_1'$ of $\calt_1$ is an extension of $\calt_0$, and
conclusion (2) holds with $\calt_1'$ in place of $\calt_1$. If
\ref{second hairy lemma}(2) holds and $n_1={n_0}$ (i.e.  $\calt_1$ is a
degenerate extension of $\calt_0$), conclusion (2) of the present
lemma holds. Finally, if \ref{second hairy lemma}(3) holds, then
$(\calt^+,\widetilde\phi^+)\in\cals$, and we have a contradiction to
the minimality of $\ADS(\widetilde\phi_1)$.
\EndProof

\Proposition\label{tower proposition}
Suppose that $g$ is an integer $\ge2$, that $M$ is a closed,
orientable $3$-manifold with $\rk M\ge g+3$, and that $\pi_1(M)$ has a
subgroup isomorphic to a genus-$g$ surface group. Then there exists a good
tower
$$\calt=(M_0,N_0,p_1,M_1,N_1,p_2,\ldots,p_n,M_n,N_n),$$ 
with base $M=M_0$, such that $N_n$ 
contains a connected incompressible closed surface of genus $\le g$.
\EndProposition

\Proof Let $K$ denote a closed, orientable surface of genus $g$. The
hypothesis implies that there is a $\pi_1$-injective PL map $\phi:K\to
M$. According to {Proposition} \ref{DS and incompressible boundary},
there exist a PL map $\widetilde\phi_0:K\to M$ homotopic to $\phi$,
and a compact, connected $3$-submanifold $N_0$ of $\inter M$, such
that (i) $\inter N_0\supset \widetilde\phi_0(K)$, (ii) the inclusion
homomorphism $\pi_1(\widetilde\phi_0(K))\to\pi_1(N_0)$ is surjective,
(iii) $\partial N_0$ is incompressible in $M$, and (iv) $N_0$ is
irreducible. According to the definitions, this means that
$\calt_0=(M,N_0)$ is a good tower of height $0$ and that
$\widetilde\phi_0$ is a good homotopy-lift of $\phi$ through
$\calt_0$.

We apply Proposition \ref{there's an extension} with these choices of
$\calt_0$ and $\widetilde\phi_0$. Conclusion (2) of \ref{there's an
extension} cannot hold since $\calt_0$ has height $0$. Hence
conclusion (1) must hold. The extension $\calt=\calt_0$ of $\calt_0$
given by conclusion (1) is a good tower whose top contains a
connected, closed incompressible surface of genus at most $g$.
\EndProof

\Lemma \label{simple tower}
Suppose that
$${\mathcal T}=(M_0,N_0,p_1,M_1,N_1,p_2,\ldots,p_n,M_n,N_n)$$ 
is a good tower of $3$-manifolds such that $M_0$ is
simple. Then the manifolds
 $M_j$ and $N_j$ are all simple for $j=0,\ldots,n$.
\EndLemma

\Proof 
By hypothesis $M_0$ is simple.  If $M_j$ is simple for a given $j\le
n$, then since $N_j$ is a submanifold of $M_j$ bounded by a (possibly
disconnected and possibly empty) incompressible surface, it is clear
from Definition \ref{simple def} that $N_j$ is simple. If $j<n$ it
then follows from \ref{simple covering} that the two-sheeted covering
space $M_{j+1}$ of $N_j$ is also simple.  \EndProof

The following theorem is the main topological result of this paper.

\Theorem\label{top 11}
Let $g$ be an integer $\ge2$.
Let $M$ be a closed simple $3$-manifold such that $\rk M \ge 4g-1$ and
$\pi_1(M)$ has a
subgroup isomorphic to a genus-$g$ surface group. Then $M$ contains a
closed, incompressible surface of genus at most
$g$.
\EndTheorem

\Proof
Applying Proposition \ref{tower proposition}, we find a good tower
$${\mathcal T}=(M_0,N_0,p_1,M_1,N_1,p_2,\ldots,p_n,M_n,N_n),$$ with
base $M_0$ homeomorphic to $M$ and with some height $n\ge0$, such that
{$N_n$} contains a connected incompressible closed surface $F$ of genus $\le
g$.  According to the definition of a good tower, $\partial N_{n}$ is
incompressible (and, {\it a priori}, possibly empty) in $M_{n}$. Hence
$N_{n}$ is $\pi_1$-injective in $ M_{n}$. Since $F$ is incompressible
in $N_{n}$, it follows that it is also incompressible in $M_{n}$.

Since $M$ is simple it follows from Lemma \ref{simple tower}
that all the $M_j$ and $N_j$ are simple.

Let $m$ denote the least integer in $\{0,\ldots,n\}$ for which $M_m$
contains a closed incompressible surface $S_m$  of genus at most
$g$.  To prove the theorem it suffices to show that $m=0$.
Let $h$ denote the genus of $S_m$.

Suppose to the contrary that $m\ge1$. Then the hypotheses of
Proposition \ref{new improved old prop 3} hold with $N_{m-1}$ and
$M_m$ playing the respective roles of $N$ and $\widetilde N$.  Suppose
that conclusion (1) of \ref{new improved old prop 3} holds, i.e. that
$ N_{m-1}$ contains an incompressible closed surface $S_{m-1}$ with
$\genus(S_{m-1})\le h\le g$.  According to the definition of a good
tower, $\partial N_{m-1}$ is an incompressible (and possibly empty)
surface in $M_{m-1}$. Hence $N_{m-1}$ is $\pi_1$-injective in $
M_{m-1}$. Since $S_{m-1}$ is incompressible in $N_{m-1}$, it follows
that it is also incompressible in $M_{m-1}$. We therefore have a
contradiction to the minimality of $m$.
  
Hence conclusion (2) of \ref{new improved old prop 3} must hold;
in particular, $N_{m-1}$ is closed, so that $N_{m-1}=M_{m-1}$; and 
$\rk M_{m-1}=\rk N_{m-1}\le4h-3\le4g-3$.
On the other hand, since
by hypothesis we have $\rk M_0\ge4g-1$, it follows from Lemma
\ref{2r-4} that for any index $j$ such that $0\le j\le n$ and such
that $M_j$ is closed, we have $\rk M_j\ge 4g-2$. This is a
contradiction, and the proof is complete.
\EndProof

\section{Proof of the geometric theorem}
\label{geometric section}

As a preliminary to the proof of Theorem \ref{geom 11} we shall point
out how the Marden tameness conjecture, recently established by Agol
\cite{agol} and by Calegari-Gabai \cite{cg}, strengthens the results
proved in \cite{accs}.

We first recall some definitions from \cite[Section 8]{accs}. Let
$\Gamma$ be a discrete torsion-free subgroup of $\Isom_+(\haitch^3)$,
and let $k\ge2$ be an integer. We say that $\lambda$ is a {\em
$k$-Margulis number} for $\Gamma$, or for $M=\haitch^3/\Gamma$, if for
any $k$ elements $\xi_1,\dots,\xi_k \in\Gamma$, and for any $z\in
\haitch^3$, we have that either
\Bullets
\item the group $\langle\xi_1,\dots,\xi_k\rangle$ is generated
by at most $k-1$ abelian subgroups, or
\item $\max_{i=1}^k \dist(\xi_i\cdot z,z)\ge \lambda$.
\EndBullets
We say that $\lambda$ is a {\em strong $k$-Margulis number} for
$\Gamma$, or for $M$, if for any $k$ elements $\xi_1,\dots,\xi_k
\in\Gamma$, and for any $z\in \haitch^3$, we have that either
\Bullets
\item the group $\langle\xi_1,\dots,\xi_k\rangle$ is generated
by at most $k-1$ abelian subgroups, or
\item $\displaystyle
\sum_{i=1}^k\frac{1}{1+e^{\dist(\xi_i\cdot z,z)}} \le
\frac{k}{1+e^{\lambda}}\ .
$
\EndBullets
We note that if $\lambda$ is a strong $k$-Margulis number for
$\Gamma$, then $\lambda$ is also a $k$-Margulis number for $\Gamma$.

A group $\Gamma$ is termed {\it $k$-free,} where $k$ is a positive
integer, if every subgroup of $\Gamma$ whose rank is at most $k$ is
free.

\Theorem\label{accs plus agol}
Let $k\ge2$ be an integer and let $\Gamma$ be a discrete subgroup of
$\Isom_+(\haitch^3)$.  Suppose that $\Gamma$ is $k$-free and purely
loxodromic.  Then $\log(2k-1)$ is a strong $k$\hyph Margulis number
for $\Gamma$.
\EndTheorem

\Proof
This is the same statement as \cite[Proposition 8.1]{accs}
except that the latter result contains the additional hypothesis that
$\Gamma$ is $k$-tame, in the sense that every subgroup of $\Gamma$
whose rank is at most $k$ is topologically tame. (To say that a
discrete torsion-free subgroup $\Delta$ of $\Isom_+(\haitch^3)$ is
topologically tame means that $\haitch^3/\Delta$ is diffeomorphic to
the interior of a compact $3$-manifold.) But the main theorem of
\cite{agol} or \cite{cg} asserts that { any finitely generated}
discrete torsion-free subgroup
$\Delta$ of $\Isom_+(\haitch^3)$ is topologically tame.
\EndProof

By combining this with another result from \cite{accs}, we shall prove:

\Theorem\label{3-free case}
Suppose that $M$ is an orientable hyperbolic $3$-manifold without
cusps and that $\pi_1(M)$ is $3$\hyph free. Then either $M$ contains a
hyperbolic ball of radius $(\log 5)/2$, or $\pi_1(M)$ is a free group
of rank $2$.
\EndTheorem

\Proof
We may write $M=\haitch^3/\Gamma$, where $\Gamma\le\Isom(\haitch^3)$
is discrete and purely loxodromic. Since $\Gamma\cong\pi_1(M)$ is
$3$-free, it follows from Theorem \ref{accs plus agol} that $\log5$ is
a strong $3$-Margulis number, and in particular a Margulis number, for
$\Gamma$ (or equivalently for $\Gamma$).  According to \cite[Theorem
9.1]{accs}, if $M$ is a hyperbolic $3$-manifold without cusps, if
$\pi_1(M)$ is $3$\hyph free and if $\lambda$ a $3$\hyph Margulis
number for $M$, then either $M$ contains a hyperbolic ball of radius
$\lambda/2$, or $\pi_1(M)$ is a free group of rank $2$. The assertion
follows.
\EndProof

\Corollary\label{3-free volume}
Suppose that $M$ is a closed orientable hyperbolic 3\hyph manifold
such that $\pi_1(M)$ is $3$-free.  Then $M$ contains a hyperbolic ball
of radius $(\log5)/2$. Hence the volume of $M$ is greater than $3.08$.
\EndCorollary

\Proof
It follows from Theorem \ref{3-free case} that either $M$ contains a
hyperbolic ball of radius $(\log5)/2$ or $\pi_1(M)$ is a free group of
rank $2$. The latter alternative is impossible, because $\Gamma$, as
the fundamental group of a closed hyperbolic $3$-manifold, must have
cohomological dimension $3$, whereas a free group has cohomological
dimension 1.  Thus $M$ must contain a hyperbolic ball of radius
$(\log5)/2$.
  
The lower bound on the volume now follows by applying B\"or\"oczky's
density estimate for hyperbolic sphere-packings as in \cite{logthree}.
\EndProof

\Theorem[Agol-Storm-Thurston]\label{firstAST}
Suppose that $M$ is a closed orientable hyperbolic $3$-manifold
containing a connected incompressible closed surface $S$. Then either
$\vol(M)>\WHAT$, or  the manifold $X$ obtained by
splitting $M$ along $S$ has the form $X=|\calw|$ for some (possibly
disconnected) book of $I$-bundles $\calw$.
\EndTheorem

\Proof
According to \cite[Corollary 2.2]{ast}, if $S$ is an
incompressible closed surface in a closed orientable hyperbolic $3$-manifold
$M$, if $X$ denotes the manifold obtained by splitting $M$ along $S$,
and if $K=\overline{X-\Sigma}$ where $\Sigma$ denotes the relative
characteristic submanifold of the simple manifold $X$, then the
volume of $M$ is greater than $\chibar(K)\cdot\WHAT$. Hence either
$\vol(M)>\WHAT$ or $\chi(K)=0$. In the latter case, we shall show that
$X$ is a book of $I$-bundles; this will complete the proof.

Note that each component of $K$ must have Euler characteristic $\le0$,
because the components of the frontier of $K$ in $X$ are essential
annuli in $X$. Since $\chi(K)=0$ it follows that each component of $K$
has Euler characteristic $0$. Hence if $Y$ denotes the union of all
components of $\Sigma$ with strictly negative Euler characteristic,
and if we set $Z=\overline{X-Y}$, then each component of $Z$ has Euler
characteristic $0$. But $Z$ is $\pi_1$-injective in $X$ because its
frontier components are essential annuli. Since $X$ is simple, the
components of $Z$ are solid tori. Since $Y=\overline{X-Z}$ is an
$I$-bundle with $Y\cap Z=\partial_vY$, and the components of
$\partial_vY$ are $\pi_1$-injective in $X$ and hence in $Z$, it
follows from the definition that $X$ is a book of $I$-bundles.
\EndProof

\Proposition\label{what if it's the other thing}
Suppose that $M$ is a closed orientable hyperbolic 3\hyph manifold
such that $\rk M\ge7$. Suppose that $\pi_1(M)$ {has a
subgroup isomorphic to a genus-$2$ surface group.} Then $\vol M\ge\WHAT$.
\EndProposition

\Proof
Since $M$ is simple and $\rk M\ge7$, it follows from Theorem
\ref{top 11} that $M$ contains either a closed incompressible surface of
genus $2$.

Suppose that  $\vol M<\WHAT$. Let $X$ denote the manifold obtained by
splitting $M$ along $S$. According to
Theorem \ref{firstAST},
 each component of $M-S$ has the form $|\calw|$ for some
book of $I$-bundles $\calw$. 

Consider the subcase in which $X$ is connected. Since $S$ has genus
$2$, we have $\chibar(X)=2$.
By Lemma \ref{moosday} it follows that
$$\rk(X)\le2\barchi(X)+1\le5.$$
Hence $$\rk M\le\rk X+1\le6,$$
a contradiction to the hypothesis.

There remains the case in which $X$ has two components, say $X_1$ and $X_2$.
Since $S$ has genus $2$,  we have
$\chibar(X_i)=1$ for $i=1,2$.
By Lemma \ref{moosday}, it follows that
$$\rk(X_i)\le2\barchi(X_i)+1=3.$$
Hence $$\rk M\le\rk X_1+\rk X_2\le6,$$
and  we have a contradiction. 
 (The bound of $6$ in this last
inequality could easily be improved to $4$, but this is obviously not
needed.) 
\EndProof

We can now prove our main geometrical result.

\Theorem\label{geom 11}
Let $M$ be a closed orientable hyperbolic $3$-manifold such that
$\vol M \le3.08$. Then $\rk M \le 6$.
\EndTheorem

\Proof Assume that $\rk M \ge 7$.  If $\pi_1(M)$ has a subgroup
isomorphic to a genus-$2$ surface group, then it follows from
Proposition \ref{what if it's the other thing} (with $g=2$)
that $\vol M\ge\WHAT>3.08$, a contradiction to the hypothesis. There
remains the possibility that $\pi_1(M)$ has no subgroup isomorphic to
a genus-$2$ surface group.  In this case, since $\rk M \ge 5$, it
follows from \cite[Proposition 7.4 and Remark 7.5]{accs} that
$\pi_1(M)$ is $3$-free. Hence by Corollary \ref{3-free volume} we have
$\vol M>3.08$, and again the hypothesis is contradicted.  \EndProof

\bibliographystyle{hplain}

\end{document}